\documentclass[11pt,reqno,intlimits,namelimits,oneside]{amsart}
\usepackage{amsfonts,amssymb,amsmath,amsthm,amssymb,amsopn,amsxtra,amstext,amsbsy}
\allowdisplaybreaks[4]
\usepackage[active]{srcltx} % Inverse Search
\renewcommand{\subjclassname}{\textup{2010} Mathematics Subject Classification}

\makeatletter
\renewcommand\subsubsection{\@startsection{subsubsection}{3}%
  \z@{.5\linespacing\@plus.7\linespacing}{-.5em}%
  {\normalfont\bfseries}}
\renewcommand\subsection{\@startsection{subsection}{2}%
  {\parindent}{.5\linespacing\@plus.7\linespacing}{-.5em}%
  {\normalfont\bfseries}}
\renewcommand\section{\@startsection{section}{1}%
  \z@{.7\linespacing\@plus\linespacing}{.5\linespacing}%
  {\normalfont\bfseries\centering}}
\renewcommand{\@secnumfont}{\bfseries}
\makeatother
\DeclareSymbolFont{AMSb}{U}{msa}{m}{n} \DeclareMathAlphabet{\msa}{U}{msa}{m}{n}
\DeclareSymbolFont{AMSb}{U}{msb}{m}{n} \DeclareMathAlphabet{\Bb}{U}{msb}{m}{n}
\DeclareSymbolFont{AMSb}{U}{eus}{m}{n} \DeclareMathAlphabet{\eusm}{U}{eus}{m}{n}
\DeclareSymbolFont{AMSb}{U}{euf}{m}{n} \DeclareMathAlphabet{\eufm}{U}{euf}{m}{n}
\DeclareSymbolFont{AMSb}{U}{eur}{m}{n} \DeclareMathAlphabet{\eurm}{U}{eur}{m}{n}
\usepackage{fancyhdr}
\newcommand{\fo}[1]{{\mbox{\footnotesize{$#1$}}}}

\numberwithin{equation}{section}
\newtheorem{theorem}{Theorem}

\newtheorem{corollary}{Corollary}
\newtheorem{proposition}{Proposition}%[section]

\newtheorem{definition}{Definition}

\newtheorem{thmx}{Theorem}

\DeclareMathOperator*{\limsp}{\overline{\lim}}

\DeclareMathOperator*{\rank}{{\rm{rank}\,}}
\DeclareMathOperator{\im}{\rm{Im}}

\usepackage{xcolor}
\usepackage{hyperref}
\hypersetup{colorlinks = false, linkbordercolor = {white}, citebordercolor = {white}}

\begin{document}
\title[Solvability of the Hankel determinant problem ... ]{{\large{
Solvability of the Hankel determinant problem for real sequences
}}}%

\author{Andrew Bakan}%
\author{Christian Berg}%
\address{Institute of Mathematics \\
National Academy of Sciences of Ukraine\\
Tereschenkivska Street~3, Kyiv 01601, Ukraine \\ }%
\email{andrew@bakan.kiev.ua}%
\address{Department of Mathematical Sciences \\ University of
Copenhagen \\ Universitetsparken 5 \\
DK-2100 Copenhagen \\  Denmark}
\email{berg@math.ku.dk}%

%\thanks{}

%44A60 Moment problems
%15A15 Determinants, permanents, other special matrix functions
%15A63 Quadratic and bilinear forms, inner products [See mainly 11Exx]
%47B36 Jacobi (tridiagonal) operators (matrices) and generalizations

%\date{^?}%
%\dedicatory{^?}%
%\commby{^?}%

\begin{abstract} To each nonzero sequence   $\eurm{s}:= \{s_{n}\}_{n \geq 0}$ of real numbers  we associate the Hankel de\-ter\-mi\-nants $D_{n} = \det \mathcal{H}_{n}$ of the Hankel matrices
$\mathcal{H}_{n}:=  (s_{i + j})_{i, j = 0}^{n}$, $n \geq 0$, and the nonempty set $\Bb{N}_{\eurm{s}}:=  \{n \geq 1 \, |\, D_{n-1} \neq 0 \}$. We also define the Hankel determinant polynomials  $P_0:=1$, and $P_n$, $n\geq 1$ as the determinant of the Hankel matrix $\mathcal H_n$ modified by replacing the last row by the monomials $1,x,\ldots,x^n$. Clearly $P_n$ is a polynomial of degree at most $n$ and of degree $n$ if and only if $n\in\Bb{N}_{\eurm{s}}$. Kronecker
established in 1881 that if $\Bb{N}_{\eurm{s}}$ is finite then $\rank \mathcal{H}_{n} = r$ for each
$n \geq r-1$, where $r := \max \Bb{N}_{\eurm{s}}$. By using an approach suggested by I.S.Iohvidov in 1969 we give a short proof of this result and a transparent proof of the
 conditions on a real sequence $\{t_n\}_{n\geq 0}$  to be of the form $t_n=D_n$, $n\geq 0$ for a real sequence $\{s_n\}_{n\geq 0}$. This is the Hankel determinant problem. We derive from the Kronecker identities   that each Hankel determinant polynomial $ P_n $  satisfying $\deg P_n = n\geq 1$
is preceded by a nonzero polynomial $P_{n-1}$ whose degree can be strictly less than $n-1$ and which has no common zeros with  $ P_n $. As an application of our results we obtain a new proof of a recent theorem by Berg and Szwarc about positive semidefiniteness of all Hankel matrices provided
that $D_0 > 0, \ldots, D_{r-1} > 0 $ and $D_n=0$ for all $n\geq r$.
\end{abstract}
\maketitle

\thispagestyle{empty}

\vspace{-0.5cm}

\section{Introduction}\label{ham}

\vspace{0.15cm}
We use the notation $\Bb N:=\{1,2,\ldots\}$ and $\Bb N_0:=\Bb N\cup \{0\}$.
To a sequence $\eurm{s}:= \{s_{n}\}_{n \geq 0}$ of real numbers  we associate  the Hankel matrices $\mathcal{H}_{n}:=  (s_{i + j})_{i, j = 0}^{n}$, $n \geq 0$ and the de\-ter\-mi\-nants $D_{n} =D_n(\eurm{s}):= \det \mathcal{H}_{n}$, $n\ge 0$. In this way we get a mapping $D: \{s_n\}_{n\ge 0}\mapsto \{D_n(\eurm{s})\}_{n\ge 0}$ in the space $\Bb R^{\Bb N_0}$ of sequences of real numbers. We call this mapping the Hankel determinant transform. It was introduced and studied by Layman in \cite{lay} who
emphasized that such a  transform  is far from being injective by proving that  a sequence  $\eurm{s}$ and its binomial transform $\beta(\eurm{s})$ defined by
$$
\beta(\eurm{s})_n:=\sum_{k=0}^n \binom{n}{k}s_k,\quad n\ge 0,
$$
have the same image under this mapping. Concerning the missing injectivity let us here just point out
that  the Hankel determinant transform of all the sequences $\{a^n\}_{n\ge 0}$, $a\in\Bb R$ is  $ \{1,0,0,\ldots\}$.

Several authors have been concerned with the sign pattern of the sequence $D(\eurm{s})$ in order to use this for the determination of the rank and signature of the Hankel matrices. This is given in rules of e.g. Jacobi, Gundelfinger and Frobenius. See \cite{gan},\cite{ioh} for a treatment of these questions, which become quite technical when zeros occur in the sequence $D(\eurm{s})$.

The Hankel determinant problem for real sequences is to characterize the image $D(\Bb R^{\Bb N_0})$ in $\Bb R^{\Bb N_0}$, i.e., to find a necessary and sufficient condition for a sequence $\eurm{t}\in\Bb R^{\Bb N_0}$ to be of the form

 \begin{gather}\label{intr1}
     \left|
\begin{array}{llllll}
s_{0} &s_{1}   &\ldots    &s_{n}   \\
s_{1}   & s_{2}  &\ldots    &s_{n+1}   \\
\ldots   & \ldots   &\ldots &\ldots  \\
s_{n}     &s_{n+1}   & \ldots   & s_{2n}
    \end{array}
\right|  = t_{n} \ , \ \ \ \ \ n \geq 0 \ .
 \end{gather}
with some sequence $\eurm{s}$ of real numbers.
 It turns out that such conditions are similar to those that were obtained by G.Frobenius  \cite[p.207]{fro}  in 1894 for all possible signs of the numbers $\{t_{n}\}_{n \geq 0}$. His arguments were simplified  by F.Gantmacher \cite[p.348]{gan} in 1959 and by I.S.~Ioh\-vi\-dov \cite[(12.8), p.83]{ioh} in 1982 in an essential way. The purpose of the present paper is to obtain a further simplification of the Frobenius reasoning by giving in Theorem~\ref{firsttheorem} a new setting
 of the approach suggested by Iohvidov in \cite[Chapter II]{ioh}. This allows to give in Section~\ref{p0th} a self-contained
proof of the following theorem.

\begin{theorem}\label{maintheorem}
   Let $\eurm{t}:= \{ t_{n} \}_{n \geq 0}$ be a   sequence of real numbers and
   \begin{gather*}
    Z_{\hspace{0.025cm} \eurm{t}} := \{\ n\geq 0  \ | \  t_{n} \neq 0 \ \} \ .
   \end{gather*}
 If $  Z_{\hspace{0.025cm} \eurm{t}} = \emptyset$ then the equation \eqref{intr1} is satisfied if and only if
 $s_{n} = 0$ for all $n \geq 0$. If $  Z_{\hspace{0.025cm} \eurm{t}} \neq \emptyset$
  consists of $1 \leq m \leq \infty$
 distinct elements $\{n_{k}\}_{0\leq k < m}$ arranged in increasing order  then the equation \eqref{intr1} is solvable if and only if the following Frobenius conditions {\rm{(}}see \cite[p.348]{gan}{\rm{)}} hold
 \begin{align*}
  & (-1)^{\tfrac{n_{0}+1}{2}} t_{n_{0}} > 0 \ , \hspace{-0.5cm} &  &  \ \ \mbox{if} \ \ \  n_{0}+1 \in 2 \Bb{N} \ ,  &
 \\
   &    (-1)^{\tfrac{n_{k+1}-n_{k}}{2}} t_{n_{k+1}}  t_{n_{k}} > 0 \ , \hspace{-0.5cm}  &  & \ \ \mbox{if} \ \ \   n_{k+1}-n_{k} \in 2 \Bb{N} \ , \ \
    0 \leq k < m -1 \ , \ \  2 \leq m \leq \infty \ .  &
 \end{align*}
 \end{theorem}

\vspace{0.3cm}
It follows from Theorem~\ref{maintheorem} that \eqref{intr1} is solvable if $t_{n} \neq 0$ for all $n \geq 0$, and not solvable if  $\eurm{t}=\{0,1,0,0,\ldots\}$. Furthermore, the condition $(-1)^{\frac{n (n+1)}{2}} t_{n}\geq 0 $  for all $n \geq 0$ is sufficient for the existence of at least one solution of \eqref{intr1}.

\smallskip
Let us formulate an elementary result about existence and uniqueness of solutions to \eqref{intr1} and which is independent of Theorem~\ref{maintheorem}. For this we need the following notation. For a  $n\times n$ determinant $A$, we denote by $A^{k,m},\;1\le k,m\le n$, the $(n-1)\times (n-1)$ determinant obtained by deleting the $k$'th row and $m$'th column of $A$. For Hankel determinants we follow
Frobenius \cite[p.212]{fro} in  writing $D_{n+1}^{\, \prime}    = D_{n+1}^{\,n+2, n+1}$, $n \geq 0$, i.e.,
\begin{gather}\label{thb4}
D_{1}^{\, \prime} = s_{1}\ , \ \
D_2^{\, \prime}  =  \left|\! \begin{array}{ll} s_{0}\!    & s_{2}
\\  s_{1} \!   &s_{3}   \end{array}\! \right| \ , \
D_{n+1}^{\, \prime}    =  \left|
\begin{array}{llllll}
s_{0} \!\! &\! s_{1} \! \! \!   &\!\!\ldots  \! & \! s_{n-2} \! & \! s_{n-1}    \! & \! s_{n+1}   \\
s_{1}  \!\!  &\! s_{2} \! \! \!  &\!\!\ldots \! & \! s_{n-1}  \! & \! s_{n}    \! & \! s_{n+2}   \\
\ldots  \!\!  &\! \ldots  \! \! \! &\!\!\ldots    \! & \! \ldots  \! & \! \ldots  \! & \! \ldots   \\
s_{n-1} \!\!     &\! s_{n} \! \! \!   &\!\! \ldots  \! & \! s_{2n-3} \! & \! s_{2n-2}   \! & \!  s_{2n} \\
s_{n}   \!\!   &\! s_{n+1} \! \! \!   &\!\! \ldots  \! & \! s_{2n-2}  \! & \! s_{2n-1}   \! & \!  s_{2n+1}
    \end{array}\!\!
\right| \  ,  \  \  n \geq 2 \  .
\end{gather}

\begin{proposition}\label{1lem} Given  two  sequences $\eurm{t}, \eurm{t}'$ of real numbers such that $t_n\neq 0$ for all $n\ge 0$, there exists a unique sequence $\eurm{s}$  of real numbers such that
  \begin{gather*}D_{n}  = t_{n} \in\Bb{R}\setminus\{0\}
 \ , \ \ \ \  D_{n+1}^{\, \prime} = t_{n}^{\, \prime}  \in\Bb{R}   \ , \ \ \ n \geq 0 \  .
  \end{gather*}
\end{proposition}

To see this we use the
Laplace expansion of  $D_{n}$ and  $D_{n+1}^{\, \prime}$  along the last  column and note that $(D'_{n+1})^{n-k,n+1}=D_n^{n-k,n+1}$. This gives the following  recurrence formulas
\begin{align*}\nonumber
 & s_{0}   =  D_{0} \ ,  \ s_{1} =   D_{1}^{\, \prime} \ ;  \hspace{0.75cm} s_{2}  D_{0} =    D_{1} + s_{1}^{2}
 \ ,
\ s_{3}  \,  D_{0} = D_{2}^{\, \prime} +  s_{1}s_{2} \ ;
\hspace{-0.2cm}  \\ &
s_{2n}  \,  D_{n-1} \hspace{0.25cm} =
  D_{n}\hspace{0.23cm} + \hspace{0.05cm}\sum\nolimits_{k=0}^{n-1} (-1)^{k}  s_{2n-1-k}  \,  D_{n}^{n-k, n+1}=: D_{n}\hspace{0.03cm} + F_n (s_{0}, ... , s_{2n-1})
\ , \hspace{-0.2cm} \\ &\nonumber
s_{2n+1} D_{n-1}\! = \!
 D_{n+1}^{\, \prime} + \hspace{0.05cm} \sum\nolimits_{k=0}^{n-1} (-1)^{k}\,  s_{2n-k}\ \ \, D_{n}^{n-k, n+1}\!  =: \! D_{n+1}^{\, \prime} +
 G_n (s_{0}, ... , s_{2n}) \ ,  \  n \!\geq\! 1 \ ;\hspace{-0.2cm}
\end{align*}
where  $D_{n}^{n-k, n+1}$, $0 \leq k \leq n-1$, depend only on
$s_{j}$, $0 \leq j \leq 2 n-1$, $F_n$ is a function of $s_0,\ldots,s_{2n-1}$ and $G_n$  a function of $s_0,\ldots,s_{2n}$.
If   $t_n'=D_{n+1}^{\, \prime}, t_n=D_{n}\neq 0,\;n\ge 0$     are assumed to be given, these relations  determine the sequence $\eurm{s}$ uniquely,
and the assertion follows.

A complete description of all solutions of \eqref{intr1}, when some of the numbers $t_{n}$ vanish, can be derived from the Frobenius  results in \cite{fro}, but this is of no relevance in the present context.

\smallskip Let $\eusm{P}[\Bb{R}]$ denote the set of  all algebraic polynomials with real coefficients.
Given a sequence $\{s_n \}_{n\geq 0}$  of real numbers  we introduce two
sequences of  polynomials in $\eusm{P}[\Bb{R}]$:

\begin{multline}\label{thb1}
 P_{0}(x)\! :=\! 1 \, , \ P_{1}(x)\! :=\! \left|\!\begin{array}{ll} s_{0}\!   & s_{1}
\\  1 \!  & x   \end{array}\!\right|
\, ,  \ \
 P_{n}(x) \!   := \!
\left|
\begin{array}{llllll}
s_{0} \!     &\!s_{1} \!   &\!s_{2}   \!    &\!\ldots \!  &\!s_{n}    \\
s_{1}  \!    &\! s_{2}  \!   &\!s_{3}   \! &\!\ldots  \!   &\! s_{n+1}   \\
\ldots  \!  &\!\! \ldots  \!      &\!\ldots \! &\!\ldots    \! &\! \ldots  \\
s_{n-1}   \! &\!s_{n}    \!  &\! s_{n+1}   \! &\!\ldots \!   &\! s_{2n-1} \\[0.15cm]
1     \!     &\! x    \!    &\! x^{2}  \!  &\!\ldots  \!   &\! x^{n}
    \end{array}
\!\!\right| \,  ; \\      Q_{0}(x) := 0 \ , \  Q_{1}(x) := s_{0}^{2} \, , \ \       Q_{n}(x)  \!  := \!
\left|
\begin{array}{llllll}
s_{0} \!\!     &\!s_{1} \!&\!s_{2}   \! \!    &\!\ldots \!    &\!s_{n}   \!    \\
s_{1}  \!\!    &\! s_{2}  \! &\!s_{3}   \! \!   &\!\ldots  \!   &\! s_{n+1} \! \\
\ldots  \!  &\!\! \ldots  \!\!     &\!\ldots \! &\!\ldots    \! &\! \ldots \! \\
s_{n-1}   \!\! &\!s_{n}    \! &\! s_{n+1}   \!\!  &\!\ldots \!   &\! s_{2n-1}\! \\[0.15cm]
0     \!\!     &\! s_{0}    \! &\!  s_{0}x\! +\! s_{1}   \!\!    &\!\ldots  \!  &\!
  \begin{textstyle}
 \sum\limits_{k=0}^{n-1} s_{k}x^{n-1-k}
  \end{textstyle}
    \end{array}\!
\right| \, , \ n\geq 1 \ .
\end{multline}

\noindent
Note that $ D_{n}^{n+1, n+1-k} =  D_{n}^{n+1-k, n+1}$, $0 \leq k \leq n$, $n \geq 1$ and
\begin{gather}\label{eq:LP}P_1(x)\! =\! D_{0} x\! -\! D_1^{\,\prime}  \, , \,
    P_n(x)\! = \! D_{n-1} x^{n}\! -\! D_n^{\,\prime} x^{n-1}\! +\!
    \sum\limits_{k=1}^{n-1} (-1)^{k-1}  x^{n-1-k}  \,  D_{n}^{n-k, n+1} \, , \,  n \geq 2 \, .
\end{gather}

 The polynomials $\{P_n\}_{n\geq 0}$ are called {\emph{Hankel determinant polynomials}} with respect to the sequence $\{s_n \}_{n\geq 0}$.
 Let $L:\eusm{P}[\Bb{R}]\to \Bb{R}$ denote the linear functional determined by
\begin{equation}\label{eq:L}
L(x^n)=s_n,\quad n\ge 0.
\end{equation}
Then
\begin{equation}\label{eq:L1}
L(x^kP_n(x))=0\quad 0\le k\le n-1,\;n\ge 1 \ ,
\end{equation}

\noindent
and also
\begin{gather}\label{thb3}
   L \left(P_n (x)^{2}\right) = D_{n} D_{n-1} \ , \ \ n \geq 0 \ , \ \ D_{-1} := 1 \ .
\end{gather}

Already Stieltjes considered this kind of functional, see \cite[p. 25]{sti2}. It is also used in \cite[Definition 2.1, p.6]{chi}).

In the classical case where all the Hankel determinants $D_n >0$, these polynomials are proportional to the classical orthonormal polynomials (see \cite[p.10; p.15; Exercise 3.1(a), p.17]{chi})
\begin{gather}\label{thb2}
    {\eurm{p}}_n (x) := \frac{P_n (x)}{\sqrt{ \vphantom{A^{A}}D_{n} D_{n-1}}} \ , \ \ n \geq 0 \ , \ \ D_{-1} := 1 \ ,
\end{gather}

\noindent
and those of the second kind.

In the general case of an arbitrary sequence $\{s_{n}\}_{n \geq 0}$ of real numbers
 Frobenius \cite[(5), p.212]{fro}  obtained  in  1894  a  recurrent relation for the  polynomials $\{P_{n}\}_{n \geq 0}$
in the following determinant form
\begin{gather}\label{tvr1}
   D_{n-1} D_{n}\ x P_n (x)\! =\!  D_{n-1}^{\,2} \,  P_{n+1}(x)\! +\! \left(D_{n-1}D_{n+1}^{\, \prime} \! -\! D_{n} D_{n}^{\, \prime} \right)\, P_{n}(x)\! +\!  D_{n}^{\,2}\, P_{n-1}(x) \ , \ \
 \end{gather}

\noindent
where $n \geq 0 $, $P_{-1} (x) \!:=\! 0 $, $P_{0}(x) \!=\! 1 $ and $D_{-1}\!:= \! 1$, $D_{0}^{\, \prime}  \!:=\! 0$.
If $D_{n} \neq 0$ for all
$n \geq 0$, then the functional $L$ is called {\emph{quasi-definite}} (see \cite[Definition 3.2, p.16]{chi}) and the monic polynomials
\begin{gather}\label{tvr7}
  p_{n} (x):= {P_{n}(x)}/{D_{n-1}} \ , \  n \geq 0 \ ,
\end{gather}

\noindent are usually considered
for which  the recurrence \eqref{tvr1} is  written in the  Jacobi form (see \cite[Theorem 4.1, p.18]{chi})
\begin{align}\label{tvr3}
& p_{n+1}(x) = \left( x -  a_{n}\right) p_{n}(x)-   b_{n} p_{n-1}(x)\  ,  \ \  n \geq 0  \ , \  \   p_{0}(x) = 1 \ , \  p_{-1} (x) = 0 \, , \\ \label{tvr4}
&  a_n =  \frac{D_{n+1}^{\, \prime}}{D_{n}}- \frac{D_{n}^{\, \prime}}{D_{n-1}}  \ , \ \  n \geq 0 \ , \ \ \
  b_{n} = \frac{D_{n}D_{n-2}}{D_{n-1}^{2}}  \ , \ \ n \geq 1 \ , \ \  a_0 = \frac{D_{1}^{\, \prime}}{D_{0}}  \ ,  \ \  b_0 = D_{0}  \ ,
\end{align}

\noindent
where the relations \eqref{tvr4} are invertible  (cp. \cite[Theorem 4.2, p.19]{chi})
\begin{gather}\label{tvr5}
    D_{n} = \prod\limits_{k=0}^{n} b_{k}^{n +1 -k} \ , \ \ D_{n+1}^{\, \prime } =\left(\sum\limits_{k=0}^{n}a_k\right) \prod\limits_{k=0}^{n} b_{k}^{n +1 -k} \ , \ \ \ \   n \geq 0 \ .
\end{gather}

\noindent
and $b_{n} \neq 0$ for all $n \geq 0$. Conversely, given the recurrence formula \eqref{tvr3} for monic
 polynomials  $\{p_n\}_{n\geq 0}$ with  two arbitrary real sequences $\{a_n\}_{n\ge 0}$ and $\{b_n\}_{n\ge 0}$ satisfying  $b_{n} \neq 0$ for all $n \geq 0$, we determine by \eqref{tvr5} and Proposition~\ref{1lem} a quasi-definite functional $L$
 such that  $L (p_n (t) p_m (t)) = 0$ and  $L (p_n (t)^{2}) \neq  0$ for all $n, m \geq 0$, $n\neq m$, by virtue of \eqref{eq:L}, \eqref{tvr7}, \eqref{eq:L1} and \eqref{thb3}. This fact is known as the
 generalized Favard theorem for quasi-definite functionals  (see \cite[Theorem 4.4, p.21]{chi}).

\smallskip
By Theorem~\ref{maintheorem}, if $\{s_n\}_{n\ge 0}$ is assumed nonzero, i.e., $s_n\neq 0$ for at least one $n\ge 0$, there exists $r\ge 1$ such that $D_{r-1}\neq 0$, and then $P_r$ is a polynomial of degree $r$. In Theorem~\ref{3th} we derive from the Kronecker identities \eqref{int11} a simple result about zeros of the  polynomials $P_r$ and $Q_r$.

\smallskip
In 1881 Kronecker \cite{kro} also characterized all those nonzero sequences   $ \{s_{n}\}_{n \geq 0}$ of real numbers  whose Hankel matrices  $( s_{i + j } )_{i, j = 0}^{\infty}$ are of finite rank, see Theorem~\ref{thb} of Section~\ref{kron}.

In Corollary~\ref{1cor} we provide a new in\-ter\-pre\-ta\-tion  of this result based on Theorem~\ref{firsttheorem}.

The results obtained by Frobenius \cite{fro} in  1894 are formulated in Theorem~\ref{thfrob} of  Section~\ref{frob}.

In Subsection~\ref{chs} we use Theorem~\ref{firsttheorem} to derive a recent theorem of Berg and Szwarc \cite{ber1}, see  Theorem~\ref{the}.
\vspace{0.25cm}

\section{Kronecker's results from 1881}\label{kron}

\vspace{0.15cm}

Let $1 \leq m \leq n$ and  $A_{n}=\{a_{i,j}\}_{i, j =1}^{n}$ be a nonzero square matrix of order $n$, where $A_n$ being nonzero means that it  has at least one nonzero element. For arbitrary
 $1 \leq i_1 < i_2 < ...< i_m \leq n$ and $1 \leq j_1 < j_2 < ...< j_m \leq n$ the determinant $\det \{a_{i_{k_{1}},j_{k_{2}}}\}_{k_{1}, k_{2} =1}^{m}$ is called a minor of $A_{n}$ of order $m$. The largest  order of the nonzero minors of $A_{n}$ is called the rank of the matrix  $A_{n}$ and is denoted by
 ${\rm{rank}}\, A_{n}$ (see {\rm{\cite[p.2]{gan}). The rank of an infinite matrix
 $A_{\infty}=\{a_{i,j}\}_{i, j =1}^{\infty}$ is defined by $ {\rm{rank}}\,A_{\infty}:=\sup\nolimits_{\,n \geq 1 } \ {\rm{rank}}\, A_{n} \in \Bb{N} \cup \{\infty\}$, where $\Bb{N}:=\{1, 2, ... \}$ (see \cite[p.205]{gan1}).

\medskip

\begin{thmx}[Kronecker (1881)]\label{thb}
Let $\{s_n \}_{n\geq 0}$ be a  nonzero sequence of real numbers, ${\mathcal{H}}_{n} :=   ( s_{i + j } )_{i, j = 0}^{n}$, $ D_{n} := \det \mathcal{H}_{n}$, $n\geq 0$, and   $\mathcal{H}_{\infty} :=   ( s_{i + j } )_{i, j = 0}^{\infty}$. Then a necessary and sufficient condition for $\mathcal{H}_{\infty}$ to have a finite rank $r\in \Bb{N}$ is that
\begin{gather}\label{tha1}
  D_{r - 1 } \neq 0 \ , \ \  \ \ D_{n} = 0 \ , \ \ \  n \geq r \ .
\end{gather}
 \end{thmx}

The necessity of the condition is formulated in Kronecker \cite[p.560]{kro}, Frobenius \cite[p. 204]{fro}, Gantmacher \cite[p.206]{gan1} and Iohvidov \cite[p. 74]{ioh}, while the sufficiency, proved by Kronecker \cite[p.563]{kro}, is less known and can be found in Iohvidov \cite[item 11, p.79]{ioh} .

\vspace{0.2cm}

It has been proved by Kronecker in  \cite[($G^{(m)}$), ($G^{\prime}$), p.567]{kro}) that
if $( s_{i + j } )_{i, j = 0}^{\infty}$ is of finite rank $r $ then
\begin{gather}\label{int5}
    {Q_{r} (x)}/{P_{r}(x) } =  \sum\nolimits_{k \geq 0 } \ { s_{k}}{x^{-k-1}}  \ \  .
\end{gather}

\noindent
Since \eqref{thb1} yields  $x^{r}P_{r}(1/x)\!\restriction_{x=0} = D_{r-1}\neq 0 $,
the change of variable $x \to 1/z$ in \eqref{int5} shows that it is equivalent to the Taylor expansion at the origin
\begin{gather*}%\label{int5a}
  \psi^{\eurm{s}}_{r} (z):=    \frac{z^{r-1}Q_{r} (1/z)}{z^{r}P_{r}(1/z) } =  \sum_{k \geq 0 }  s_{k} z^{k}
\end{gather*}

\noindent
of the analytic function $ \psi^{\eurm{s}}_{r}$ on the open disk $|z|< 1/\rho_{r}$ where $\rho_{r} := \max \{\ |z| \ | \ P_{r}(z) = 0\}$. Therefore the series in the righthand side of \eqref{int5} converges absolutely for every $|x| > \rho_{r}$, and
\eqref{int7} below holds by the Cauchy-Hadamard formula (see \cite[(2), p.200]{rud}).

Conversely, Kronecker proved in  \cite[p.568]{kro} that if the numbers $\{s_n \}_{n\geq 0}$ are the coefficients in the expansion \eqref{int8}  of $q/p$ for   $p,  q \in \eusm{P}[\Bb{R}]$, $\deg p = r \in \Bb{N}$ and $\deg q < r$ then $\mathcal{H}_{\infty}$ has the  rank $r$, provided  $D_{r-1}\!\neq\! 0$ (see  \cite[Section 45, p.198]{pra}, \cite[Theorem 8, p.207]{gan1}). Thus, the following characterization of the Hankel matrices of finite rank holds.

\vspace{0.2cm}
\begin{thmx}\label{thc}
Let $\{s_n \}_{n\geq 0}$ be a  nonzero sequence of real numbers and    $\mathcal{H}_{\infty} :=   ( s_{i + j } )_{i, j = 0}^{\infty}$.

\vspace{0.25cm}\label{thca}
{\rm{{\bf{(a)}}}} If $\mathcal{H}_{\infty}$ has a finite rank $r\in \Bb{N}$ then $\deg P_{r} = r$,
\begin{gather}\label{int7}
    \limsp_{k\to \infty}\sqrt[k]{\left|s_{k} \right|} =   \max \{\ |z| \ | \ z \in \Bb{C} \ , \  P_{r}(z) = 0\ \}   \ ,
\end{gather}

\noindent
and  \vspace{-0.3cm}
\begin{gather}\label{int6}
 \sum_{k \geq 0 } \frac{ s_{k}}{x^{k+1}}   = \frac{Q_{r} (x)}{P_{r} (x) } \ ,
\end{gather}

\noindent
where the series is absolutely convergent for every $|x| >  \max \{\ |z| \ | \ z \in \Bb{C} \ , \  P_{r}(z) = 0\ \} $.

\vspace{0.25cm}\label{thcb}
{\rm{{\bf{(b)}}}} If $R := \limsp_{k\to \infty}\sqrt[k]{\left|s_{k} \right|} < +\infty $ and there exist
$p,q \in\eusm{P}[\Bb{R}]$, $p$ of degree $r \in \Bb{N}$ and $q$ of degree at most $r-1$ such that
\begin{gather}\label{int8}
   \sum_{k \geq 0 } \frac{ s_{k}}{z^{k+1}}   =     \frac{q (z)}{p (z) } \ , \ \ \ \ |z| > R \ ,
\end{gather}

\noindent
then $\rank \mathcal{H}_{\infty} \leq r$, where the equality is attained if $p$  and $q$ have no common roots.
 \end{thmx}

\vspace{0.2cm}  The following    theorem of Kronecker \cite[pp.560, 561, 571]{kro} clarifies the  structure of the sequences satisfying  ${\rm{rank}}\,\mathcal{H}_{\infty} < \infty $ (see also \cite[Theorem 7, p.205]{gan1} and \cite[p.234]{gan1}).

\vspace{0.2cm}
\begin{thmx}\label{thd}
Let $\{s_n \}_{n\geq 0}$ be a  nonzero sequence of real numbers and
$\mathcal{H}_{\infty} :=   ( s_{i + j } )_{i, j = 0}^{\infty}$.

\vspace{0.25cm}
{\rm{{\bf{(a)}}}} $\mathcal{H}_{\infty}$ has a finite rank $r\in \Bb{N}$ if and only if $D_{r-1}\neq 0$ and there exist $r$ numbers \\  $d_{0}$,
 $d_{1}$, ... , $d_{r-1}$ such that
\begin{gather}\label{int9}
\sum\limits_{k=0}^{r-1}\ d_{k}\ s_{k + m} = \ s_{r  +  m} \ ,  \ \ \ \ \  m \geq 0 \ .
\end{gather}

\vspace{0.25cm}
{\rm{{\bf{(b)}}}} If $\mathcal{H}_{\infty}$ has a finite rank $r\!\in\! \Bb{N}$ then for every $n\geq 0$  there exist $r$ numbers  \\ $d_{n, 0}$,  $d_{n, 1}$, ..., $d_{n, r-1}$ such that
\begin{gather*}%\label{int10}
\sum\limits_{k=0}^{r-1}\ d_{n, k}\ s_{k + m} = \ s_{r + n +  m} \ ,  \ \ \ \ \  m \geq 0 \ ,
\end{gather*}

\noindent
where $d_{0, k}$ is equal to $d_{k}$ from \eqref{int9}  for each  $0 \leq k \leq r-1$.

\vspace{0.3cm}
{\rm{{\bf{(c)}}}}  If $\mathcal{H}_{\infty}$ has a finite rank $r\!\in\! \Bb{N}$ then  the sequence $\{s_n \}_{n\geq 0}$ is uniquely determined by the values of  $s_{0}$,  $s_{1}$, ..., $s_{2r-1}$.
\end{thmx}

Finally, we note that the equality \eqref{int6} proved by Kronecker in
\cite[($G^{(m)}$), ($G^{\prime}$), p.567]{kro}) asserts implicitly that  the polynomials $P_{r}$ and $Q_{r}$ have no common roots  provided that  $D_{r-1}\neq 0$. Furthermore, this fact also follows  from the identity
 \begin{gather}\label{int11}
  P_{r-1} (x) Q_{r} (x) -  P_{r} (x) Q_{r-1} (x) = D_{r-1}^{2} \ ,
 \end{gather}

\noindent
written by  Kronecker in \cite[(F), p.564]{kro} for arbitrary $ r \geq 1$ (see also \cite[(14), p.220]{fro}).

Observe that \eqref{int11} can easily be proved when $D_{n}\neq 0$ for all $n \geq 0 $ (see \cite[III.15, p.48]{ger}, \cite[Theorem 2.12, p.54]{bre}).  These restrictions can  be removed by the so-called perturbation
technique. More precisely, by using the  Hilbert matrix $\eusm{M}_{n}^{1} := ((i + j +1)^{-1})_{i, j =0 }^{n}$, for every
$\varepsilon > 0$ we introduce the perturbed sequence
\begin{gather*}%\label{int11a}
 \{s_{n}^{\varepsilon} \}_{n \geq 0} \ , \ \ \ \
    s_{n}^{\varepsilon} := s_{n} + \frac{\varepsilon^{n + 1}}{n + 1} \ , \  \eusm{M}_{n}^{\varepsilon}:= \Big(\,
    \frac{\varepsilon^{i + j +1}}{i + j +1}\,\Big)_{i, j =0 }^{n} \ , \ \ n   \geq 0  \ ,
\end{gather*}
whose Hankel  determinant
$D_{n}^{\varepsilon} = \det ( {\mathcal{H}}_{n} +  \eusm{M}_{n}^{\varepsilon} ) =
{\eurm{m}}_n \varepsilon^{(n+1)^{2}} + ... $ for every $n\geq 0$ is a polynomial  of degree $(n+\!1)^{2}$ in the variable $\varepsilon$ with positive  leading coefficient ${\eurm{m}}_n := \det  \eusm{M}_{n}^{1} > 0$ (see \cite[3, p.92]{psz1}) . Since the zeros of all polynomials  $D_n^{\varepsilon}$, $n \geq 0$, form an at most countable set, there exists a sequence $\{\varepsilon_{k}\}_{k\geq 0}$ of positive numbers $\varepsilon_{k}$ tending to zero as $k \to \infty $ such that
\begin{gather*}%\label{int11b}
\det \big\{s_{i+j}^{\,\varepsilon_{k}}\big\}_{i, j =0}^{n}
\neq 0  \ , \ \  n, k \geq 0 \ .
 \end{gather*}

\noindent
With \eqref{int11} in hand for $\{s_{n}^{\,\varepsilon_{k}} \}_{n \geq 0}$, $k \geq 0$, we conclude by the  continuous dependence of  \eqref{int11} on $s_{m}^{\,\varepsilon_{k}} \ $, $0 \leq m \leq 2r-1$, that  \eqref{int11} holds for  $\{s_{n}\}_{n \geq 0}$.

\smallskip
It also follows from \eqref{int11} that  $P_{r}$ and $ P_{r-1}$ have no common roots, provided  that  $D_{r-1}\neq 0$. We have therefore proved the following property (cf. \cite[Theorem 2.14, p.57]{bre}).

\begin{theorem}\label{3th}
Let $\{s_{n}\}_{n \geq 0}$ be an arbitrary nonzero sequence of real numbers and $r$ be a positive integer satisfying $D_{r-1}\neq 0$. Then $\deg P_{r} = r$, $P_{r-1} \not\equiv 0$, $Q_{r} \not\equiv 0$ and
 the polynomial $P_{r}$  has no common zeros with the polynomials $P_{r-1}$ and $ Q_{r}$.
\end{theorem}

\medskip
Observe, that Theorem~\ref{3th} can also be easily deduced from \cite[Theorem 1.9, p.80; Theorem 1.3(ii), p.44]{dra}.
 We will in the sequel use the following notion.

\begin{definition} {{Let $\{s_n \}_{n\geq 0}$ be a  nonzero sequence of real numbers. The rank
  of the infinite Hankel matrix $ ( s_{i + j } )_{i, j = 0}^{\infty}$ is called the {{Hankel rank}} of $\{s_n \}_{n\geq 0}$. }}
\end{definition}

\noindent
Since ${\rm{rank}}\, ( s_{i + j } )_{i, j = 0}^{\infty} \in \Bb{N} \cup \{\infty\}$, the Hankel rank of a real nonzero sequence can be equal to any positive integer or infinity.

\vspace{0.25cm}

\section{Frobenius' theorem from 1894}\label{frob}

\vspace{0.15cm}

 Let $\eurm{s}:= \{s_n \}_{n\geq 0}$ be an arbitrary nonzero sequence of real numbers
and
\begin{gather}\label{ioh1}
    \Bb{N}_{\eurm{s}} := \left\{\  r \in \Bb{N} \ \big| \ D_{r-1} \neq 0 \ \right\} \ .
\end{gather}

\noindent  Theorem~\ref{maintheorem} yields $ \Bb{N}_{\fo{\eurm{s}}} \neq \emptyset $.
Suppose  that $ \Bb{N}_{\fo{\eurm{s}}}$  consists of $m$ ($1 \leq m \leq \infty$)
 distinct elements $\{n_{k}\}_{1\leq k < m+1}$ arranged in increasing order  and  $n_0 := 0$, where it is assumed that $a+ \infty = \infty$ for arbitrary $a\in \Bb{R}$.
  Then
 \begin{gather}\label{ioh1a}
 \{0\} \cup  \Bb{N}_{\eurm{s}} = \{n_{k}\}_{0\leq k < m+1} \ , \ \ 1 \leq m \leq \infty \ , \ \
 0 = n_0 < n_{1} < \ldots \ \ .
 \end{gather}

We say that the  Hankel determinant
polynomial $P_n$ defined by \eqref{thb1} is of {\emph{full degree}} if $\deg P_n = n$.
It follows from \eqref{thb1}, \eqref{eq:LP}, \eqref{ioh1} and  \eqref{ioh1a}
that $P_n$ is of full degree if  and only if $n= n_{k}$ for some $0 \leq k < m+1$, i.e.,
\begin{gather}\nonumber
    \left\{ \ P_n \ | \ \deg P_{n} = n \ , \ n \geq 0 \ \right\} = \left\{P_{n_{k}}\right\}_{0\leq k < m+1} =
    \left\{P_{n_{0}}\equiv 1 \, , \, P_{n_{1}} \, , \, \ldots  \, , \, P_{n_{k}}  \, , \, P_{n_{k+1}}  \, , \, \ldots \ \right\} \ ,
    \\[0.2cm]
     \deg P_{n_{k}} = n_{k} \ , \  0 \leq k < m+1\ .
\label{ioh1b}\end{gather}

Theorem~\ref{3th} states that the identities \eqref{int11} proved by Kronecker in 1881 imply  that for each  $0 \leq k < m$   the polynomial $P_{n_{k+1}}$
 is preceded by a nonzero polynomial $P_{n_{k+1}-1}$ which has no common zeros with  $P_{n_{k+1}}$ and whose degree can be strictly less than $n_{k+1}-1$.

 In 1894 Frobenius established  \cite[(10), p.210]{fro}  that  $P_{n_{k+1}-1}$ for such $k$ is proportional with a nonzero real constant of proportionality  to the previous polynomial  $P_{n_{k}}$ of full degree provided that
$\deg P_{n_{k+1}-1}< n_{k+1}-1$, i.e., there exists $\gamma_{k} \in \Bb{R}\setminus \{0\}$ such that
\begin{gather*}
  P_{n_{k+1}-1}(x) = \gamma_{k}P_{n_{k}} (x) \ ,
\end{gather*}

\noindent
if $n_{k+1} - n_{k} \geq 2$ and $0 \leq k < m$ (see also  \cite[Theorem 1.3(ii), p.44]{dra}). Furthermore, he proved in \cite[(8), p.214]{fro} that for $m \geq 2$ the recurrence relations
\begin{gather*}
   p_{n_{k+1} }(x) = a_{ k} (x) p_{n_{k} }(x) - \beta_{k} p_{n_{k-1} }(x) \ , \ 1 \leq k < m \ , \
\end{gather*}

\noindent  hold  between the monic polynomials
\begin{gather*}
 p_{n_{k} }(x) := P_{n_{k}}(x)/ D_{n_{k}-1} \ , \ 0\leq k < m+1 \ , \  \ \ \ D_{-1}:=  1\ , \
\end{gather*}

\noindent corresponding to the polynomials $P_n$ of full degree,
 where
 $\{\beta_{k}\}_{1\leq k < m}$ are  nonzero real numbers  and $a_{ k} (x)\in  \eusm{P}[\Bb{R}]$ is a monic polynomial of degree $n_{k+1} - n_{k} $ for every $1\leq k < m $ (see also \cite[Remark 1.2, p.71]{dra}).  It is also proved in
  \cite[(9), p.210]{fro} that
 \begin{gather*}
  P_{n_{k}+1}(x) \equiv   ...  \equiv  P_{n_{k+1}-2}(x)\equiv 0
 \end{gather*}

\noindent   provided that
 $n_{k+1} - n_{k} \geq 3$ and $0\leq k < m $ (see also \cite[Theorem 1.3, p.44]{dra}). Thus,  the following  theorem was proved by Frobenius \cite{fro}  in 1894 .

 \vspace{0.15cm}
\begin{thmx}\label{thfrob}
  Let  $\{s_n \}_{n\geq 0}$ be an arbitrary nonzero sequence of real numbers and $\Bb{N}_{\eurm{s}} $, $ m$ and $\{n_{k}\}_{0\leq k < m+1}$ be defined as in \eqref{ioh1} and \eqref{ioh1a}.

 \smallskip
  For  the Hankel determinant polynomials  $\{P_{n}\}_{n \geq  0}$ defined by \eqref{thb1} the following assertions hold.

\vspace{0.35cm}\noindent
{\bf{(a)}} If $n_{1}\geq 2$ then   there exists $\gamma_{0}\in \Bb{R}\setminus \{0\}$ such that

\vspace{-0.25cm}
 \begin{gather*}
 P_{0} \equiv 1\ , \    \ P_{1}  \equiv \gamma_{0} \ , \ \ \deg P_{n_{1}} =  n_{1} = 2 \ ,
\end{gather*}

 \vspace{0.15cm}
\noindent
when $n_{1}= 2$ and

\vspace{-0.25cm}
   \begin{gather*}
  P_{0} \equiv 1\ , \    P_{1} \equiv 0 \ , \ \ \ldots \ ,\ P_{n_{1} -2}  \equiv 0 \ , \ \ P_{n_{1} -1}  \equiv \gamma_{0} \ , \ \ \deg P_{n_{1}} = n_{1} \ ,
\end{gather*}

 \vspace{0.15cm}
\noindent
when $n_{1}\geq 3$.

\vspace{0.45cm}\noindent
{\bf{(b)}}  If $m \geq 2$,  $1\leq k < m$ and $n_{k+1} - n_{k}\geq 2$ then there exists $\gamma_{k}\neq 0$ such that

\vspace{-0.25cm}
\begin{gather*}
  \deg P_{n_{k}} = n_{k} \ , \ \  P_{n_{k+1} -1} = \gamma_{k}  P_{n_{k}} \ , \ \ \deg P_{n_{k+1}} = n_{k+1} \ ,
\end{gather*}

 \vspace{0.2cm}
\noindent
when $n_{k+1} - n_{k}= 2$ and

\vspace{-0.25cm}
\begin{gather*}
  \deg P_{n_{k}} = n_{k} \ , \ \ P_{n_{k} +1} \equiv 0 \ , \ \ ... \ \  ,   \   P_{n_{k+1} -2}  \equiv 0 \ , \ \ P_{n_{k+1} -1} = \gamma_{k}  P_{n_{k}} \ , \ \ \deg P_{n_{k+1}} = n_{k+1} \ ,
\end{gather*}

 \vspace{0.2cm}
\noindent
when $n_{k+1} - n_{k}\geq 3$.

\vspace{0.55cm}\noindent
{\bf{(c)}}  If $m \geq 2$ then for the  monic polynomials
\begin{gather*}
    p_{0} (x) = 1 \ , \   \  p_{n_{k}} (x) := \frac{P_{n_{k}} (x)}{D_{n_{k}-1}} \ , \ \ 0\leq k < m+1 \ ,
\end{gather*}

 \vspace{0.15cm}
\noindent
there exist  monic polynomials  in $\eusm{P}[\Bb{R}]$

\vspace{-0.25cm}
\begin{gather*}
    a_{k} (x)  \ ,  \ \
     \deg  a_{k} (x) = n_{k+1}-n_{k} \geq 1 \ ,  \  \ \ 0 \leq k < m \ ,
\end{gather*}

 \vspace{0.15cm}
\noindent
and  nonzero real numbers  $\{\beta_{k}\}_{0\leq k < m}$ such that

\vspace{-0.25cm}
\begin{gather}\label{1fthfrob}
 p_{n_{k+1}} (x) = a_{k} (x) p_{n_{k}} (x) - \beta_{k} p_{n_{k-1}} (x) \ , \ \ 0\leq k < m \ , \  \  \  \ p_{n_{-1}} := 0  \ .
\end{gather}

\vspace{0.55cm}\noindent
{\bf{(d)}}  If $m < \infty$ then $ n_{m} = \max \Bb{N}_{\eurm{s}}$ and $P_{n} \equiv 0$ for all $n \geq n_{m} + 1$.

\end{thmx}

\vspace{0.3cm}
It should  be noted that Theorem~\ref{thfrob}\! (d)  follows directly from
Theorem~\ref{thb} and Theorem~\ref{thd}\!~(a). Indeed,   the conditions of
Theorem~\ref{thfrob}\! (d)  imply the validity of \eqref{tha1}
 for $r=n_{m}$ and in view of Theorem~\ref{thb} we obtain that $\mathcal{H}_{\infty}$ has a finite rank $n_{m}$.
But for arbitrary $n \geq n_{m} + 1$  the $(n_m +1)$-th row of the determinant for $P_n$  in \eqref{thb1}  is the linear combination of the first $n_m$ rows by virtue of \eqref{int9}. Hence, $P_{n} \equiv 0$ and
 the desired result is proved.

 \smallskip Theorem~\ref{thfrob} shows that
  except of polynomials of full degree and proportional to them the sequence $\{P_n\}_{n \geq 0}$ defined in
   \eqref{thb1} contains no other  nonzero polynomials. Furthermore, if $n \geq 1 $ then it follows from  $P_n \equiv 0$ that $\deg P_{n+1} < n+1$ while $P_n \not\equiv 0$ and $\deg P_{n} < n$ imply $\deg P_{n+1} = n+1$. Observe that  Theorem~\ref{thfrob}\! (c) was essentially  generalized by A.Draux \cite[Theorem 6.2, p.477]{dra} in 1983.

\vspace{0.25cm}

\section{Iohvidov's approach  from 1969}

\vspace{0.15cm}
Throughout this section we fix an arbitrary nonzero sequence   $\eurm{s}:= \{s_{n}\}_{n \geq 0}$  of real numbers
and use the set $ \Bb{N}_{\fo{\eurm{s}}}$ defined  as in \eqref{ioh1}.
The analysis below will not use the statements from the previous Sections~\ref{kron} and~\ref{frob}.

\smallskip
 In 1969 Iohvidov \cite{ioh1} (see also \cite{ioh}) suggested a new technique for dealing with  Hankel matrices.
For every $r \in  \Bb{N}_{\fo{\eurm{s}}} $ he proposed to use the approximating
sequence $\eurm{s}^{(r)}$  defined as follows.

\medskip
 We first  put
\begin{gather}\label{ioh2}
  s_{n}^{(r)} = s_{n} \ , \ 0 \leq n \leq 2 r -1 \ .
\end{gather}

\noindent
Since the first $2 \, r-1$ numbers  $s_{0}$, $s_{1}$, $s_{2}$, $\ldots$ $s_{2 r-2}$ of the sequence $\eurm{s}$ satisfy
\begin{gather}\label{ioh3}
  D_{r-1} =  \begin{vmatrix} s_{0} &  s_{1}& \ldots &  s_{r-2} &  s_{r-1} \\  s_{1} & s_{2}& \ldots &  s_{r-1}&  s_{r} \\     \ldots  &  \ldots &  \ldots &  \ldots &  \ldots  \\
s_{r-2} & s_{r-1}& \ldots &  s_{2r-4}&  s_{2r-3}  \\    s_{r-1} & s_{r}& \ldots &  s_{2r-3}&  s_{2r-2}  \\ \end{vmatrix}   \neq 0 \ ,
\end{gather}

\noindent
it is possible to determine uniquely all $r$ numbers $d_{0}^{(r)}$, $d_{1}^{(r)}$,  ... , $d_{r-1}^{(r)}$ from the system
\begin{gather}\label{ioh4}
   \begin{pmatrix} s_{0} &  s_{1}& \ldots &  s_{r-2} &  s_{r-1} \\[0.1cm]
     s_{1} & s_{2}& \ldots &  s_{r-1}&  s_{r} \\
          \ldots  &  \ldots &  \ldots &  \ldots &  \ldots  \\[0.1cm]
          s_{r-2} & s_{r-1}& \ldots &  s_{2r-4}&  s_{2r-3}  \\[0.1cm]
              s_{r-1} & s_{r}& \ldots &  s_{2r-3}&  s_{2r-2}
              \\ \end{pmatrix}
\begin{pmatrix}
 d_{0}^{(r)} \\  d_{1}^{(r)} \\
  \ldots \\ d_{r-2}^{(r)} \\  d_{r-1}^{(r)} \\
\end{pmatrix} = \begin{pmatrix}
 s_{r} \\[0.1cm]  s_{r+1} \\
  \ldots \\[0.1cm] s_{2r-2} \\[0.1cm]  s_{2r-1} \\
\end{pmatrix} \ ,
\end{gather}

\noindent
and then to define recursively the numbers $s_{2 r}^{(r)}$, $s_{2 r+1}^{(r)}$, ... by the formulas
\begin{gather}\label{ioh5}
s_{2 r  }^{(r)}\! =\! s_{r} d_{0}^{(r)}\! +\! s_{r+1} d_{1}^{(r)}\! +\! \ldots\! +\!
s_{2r-2} d_{r-2}^{(r)}\! +\! s_{2r-1} d_{r-1}^{(r)} \, , \
   s_{2 r + m }^{(r)} \!= \!\sum\limits_{k=0}^{r-1} \ s_{r + m + k}^{(r)}  d_{k}^{(r)} \, ,  \    m \geq 1  \, . \end{gather}

\noindent
We obtain the sequence
\begin{gather}\label{ioh6}
 {\eurm{s}}^{(r)}:=  \{s_{n}^{(r)}\}_{n \geq 0} \ , \ \ \ \ {\eurm{s}}^{(r)} =\{ s_{0} \ , \ s_{1}\ ,  \ s_{2} \ , \ \ldots \ , \ s_{2 r-1} \ , \
s_{2 r}^{(r)}\ , \ s_{2 r+1}^{(r)} \ , \  ... \  \}  \ ,
\end{gather}
whose terms satisfy  a homogeneous linear recurrence relation with constant coefficients of the form
\begin{gather}\label{ioh7}
 \sum\nolimits_{k=0}^{r-1}\ d_{k}^{(r)}  s_{k + m}^{(r)} =    s_{r + m}^{(r)}\ , \ \   m \geq 0  \ ,
\end{gather}

\noindent
which is equivalent to the simultaneous validity of \eqref{ioh5} and  \eqref{ioh4} provided that \eqref{ioh2} holds.

\bigskip
The statements of Theorem~\ref{firsttheorem}\! (d) and (e) below   follow easily from  the main results of Iohvidov in Chapter II of \cite{ioh}, while their proofs given in Subsections~\ref{p1thd} and~\ref{p1the} are based on  a somewhat different approach than the one used in \cite{ioh}.

\vspace{0.15cm}
\begin{theorem}\label{firsttheorem}
Let $\eurm{s}\!:= \!\{s_n \}_{n\geq 0}$ be a  nonzero sequence of real numbers and  ${\mathcal{H}}_{n} :=   ( s_{i + j } )_{i, j = 0}^{n}$, $ D_{n} := \det \mathcal{H}_{n}$, $n\geq 0$.

\vspace{0.35cm}
{\rm{{\bf{(a)}}}} \ \ For arbitrary  $n\in\Bb N$ the property

\vspace{-0.25cm}
\begin{gather}\label{1th1}
    s_{0}= s_{1} = s_{2} =\ \ldots \  = s_{n -1} = 0  \ , \ \  s_{n} \neq 0 \ ,
\end{gather}

\noindent
is equivalent to
\begin{gather*}%\label{1th2}
    D_{0}= D_{1}  =\ \ldots \  = D_{n -1} = 0  \ , \ D_{n}\neq 0   \ .
   \end{gather*}

\vspace{0.4cm}
\noindent
If $n \in \Bb{N}$ and  $ D_{0}= D_{1}  =\ \ldots \  = D_{n -1} = 0 $ then

\begin{gather*}%\label{1th2a}
    D_{n} = (-1)^{\frac{n (n+1)}{2}} s_{n}^{n+1} \ \ \ .
\end{gather*}

\vspace{0.45cm}
{\rm{{\bf{(b)}}}} \ \
The set \vspace{-0.2cm}
\begin{gather}\label{1th3}
   \Bb{N}_{\fo{\eurm{s}}} := \left\{\  r \geq 1 \ \big| \ D_{r-1} \neq 0 \ \right\} \
\end{gather}

 \vspace{0.15cm}
\noindent
is nonempty and for every $r \in \Bb{N}_{\fo{\eurm{s}}}$ the formulas \eqref{ioh2}, \eqref{ioh4} and \eqref{ioh5}
produce the sequence

\vspace{-0.25cm}
\begin{gather*}%\label{1th4}
   {\eurm{s}}^{(r)}:=  \{s_{n}^{(r)}\}_{n \geq 0} \ , \ \ \ \ {\eurm{s}}^{(r)} =\{ s_{0} \ , \ s_{1}\ ,  \ s_{2} \ , \ \ldots \ , \ s_{2 r-1} \ , \
s_{2 r}^{(r)}\ , \ s_{2 r+1}^{(r)} \ , \  ... \  \}  \ ,
\end{gather*}

 \vspace{0.15cm}
\noindent
such that\vspace{-0.2cm}
\begin{gather}\label{1th5}
 {\rm{rank}}\,  \big( \, s_{i + j }^{(r)}\, \big)_{i, j = 0}^{\infty} = r \ .
\end{gather}

\vspace{0.35cm}
{\rm{{\bf{(c)}}}} \ \ For every $r \in \Bb{N}_{\fo{\eurm{s}}}$ we have
\begin{gather}\label{1th6}
    \limsp_{k\to \infty}\sqrt[k]{\left|s_{k}^{(r)} \right|} =
     \max \{\ |z| \ | \ z \in \Bb{C} \ , \  P_{r}(z) = 0\ \}  \ ,
\end{gather}

\noindent
and
\begin{gather}\label{1th7}
  \frac{s_{0}}{x} +  \frac{s_{1}}{x^{2}} + ... + \frac{s_{2r-1}}{x^{2r}} +
    \sum_{k \geq 2r } \frac{ s_{k}^{(r)}}{x^{k+1}}  = \frac{Q_{r} (x)}{P_{r} (x) } \ ,
\end{gather}

\vspace{0.15cm}
\noindent
where the series is absolutely convergent for every $|x| > \max \{\ |z| \ | \ z \in \Bb{C} \ , \  P_{r}(z) = 0\ \} $.

\vspace{0.5cm} 
{\rm{{\bf{(d)}}}} \ \ For arbitrary  $d \in \Bb{N}$  and $r \in \Bb{N}_{\fo{\eurm{s}}}$ the following statements hold:

\vspace{-0.25cm}
\begin{align}\label{1th8}
 & & D_{r }\!=\!  \left(s_{2r} - s_{2r}^{(r)}\right)   D_{r-1} \ ; \!\!\!  & & & & &\!\!\!
 D_{r+1}^{\, \prime}    \!=\!   \left(s_{2r+1}-s^{(r)}_{2r+1}\right)D_{r-1} -  \left( s_{2r} -s^{(r)}_{2r}\right)D_{r}^{\, \prime} \ ;
 \\[0.2cm] \label{1th9}& &  D_{r}\! =\! ...\!=\! D_{r + d-1}\!=\!0\!\!\!  & & \Leftrightarrow & & & \!\!\!
s_{2r}\! =\!  s_{2r}^{(r)} \, , \  s_{2r+1}\! =\!  s_{2r+1}^{(r)} \, , \ \ldots \ s_{2r+d-1}\! = \! s_{2r+d-1}^{(r)}
  \ ; \\[0.3cm] \label{1th10}& &  D_{r}\! =\! ...\!=\! D_{r + d-1}\!=\!0\!\!\!  & & \Rightarrow & & &\!\!\!  D_{r +  d} \!= \!(-1)^{{\tfrac{d(d  +1)}{2}} } \left(s_{2r+d} - s_{2r+d}^{(r)}\right)^{d  +1}   D_{r-1}
  \ .\end{align}

\vspace{0.5cm}
{\rm{{\bf{(e)}}}} \ \ For every $r \in \Bb{N}_{\fo{\eurm{s}}}$  the Iohvidov characteristic function

\vspace{-0.3cm}
\begin{gather}\label{1th11}
   d_{r} \ := \ \inf \left\{ \ m \geq 0 \ \ \big|\ \ s_{2 r  +m} \neq {s}_{2 r +m}^{(r)} \  \right\} \ \in \ \{ 0, 1, 2, .... \} \cup \left\{\infty \right\} \ ,
\end{gather}

\vspace{0.15cm}
\noindent where it is assumed that $\inf \emptyset := \infty $,
 possesses the following property

 \vspace{-0.3cm}
 \begin{gather}\label{1th12}
   d_{r}\ =\ \inf \left\{\ \ m \geq 0 \ \  \big| \ \ D_{ r  +m} \neq  0\ \   \right\}   \ .
 \end{gather}

\noindent
In particular, for arbitrary $r \in \Bb{N}_{\fo{\eurm{s}}}$ we have

 \vspace{-0.3cm}
 \begin{gather}\label{1th13}
   d_{r} = \infty \  \  \Leftrightarrow \ \ \eurm{s} =  {\eurm{s}}^{(r)}\ \ \Leftrightarrow \ \
D_{r-1}\neq 0 \ , \ \ D_{r} = D_{r + 1} = ...= 0   \ .\end{gather}
 \end{theorem}

 \vspace{0.2cm}
 The case of Theorem~\ref{firsttheorem}(a) was first considered by Frobenius \cite[p.206]{fro} in 1894.
 In 1969 Iohvidov \cite[(5), p.244]{ioh1}  introduced   the characteristic function $d_{r}$ and established in \cite[(7), p.246]{ioh1}  the equalities \eqref{1th10}  (see also
 \cite[(10.5), p.62; (11.2), p.70]{ioh}). A formula similar to \eqref{1th10} has been established in \cite[Lemma 2.3]{ber1}. However,
the setting of Theorem~\ref{firsttheorem}\!  (d)  and (e) differs from that of
 \cite[Chapter 2]{ioh} because only the finite Hankel matrices $( s_{i + j } )_{i, j = 0}^{n}$ are   considered there.

\vspace{0.25cm}

\section{\texorpdfstring{Consequences of Theorem~\ref{firsttheorem}}{Consequences of Theorem 3} }

 \vspace{0.2cm}

 \subsection{\texorpdfstring{Approximating sequence}{Approximating sequence}}
 The first immediate consequence of Theorem~\ref{firsttheorem} is that the sequence $ {\eurm{s}}^{(r)}$ for every
 $r \in \Bb{N}_{\fo{\eurm{s}}}$ can equivalently be defined by the expansion \eqref{1th7} which in view of
 $x^{r}P_{r}(1/x)\!\restriction_{x=0} = D_{r-1}\neq 0 $  can be considered as the Taylor expansion at the origin of the rational function
 \begin{gather*}%\label{ioh10}
   \frac{z^{r-1}Q_{r} (1/z)}{z^{r}P_{r}(1/z) } =  \sum\nolimits_{k \geq 0 } \  s_{k}^{(r)} z^{k} \ , \ \ |z| < \min \{\ |\zeta| \ | \ \zeta \in \Bb{C}\setminus\{0\} \ , \ \  P_{r}(1/\zeta) = 0\ \} \ .
 \end{gather*}

\medskip
 \subsection{\texorpdfstring{Finiteness of rank}{Finiteness of rank}} Assume now that for a given nonzero sequence $\{s_n \}_{n\geq 0}$ of real numbers  the infinite Hankel matrix $ ( s_{i + j } )_{i, j = 0}^{\infty}$ has a finite rank. Then the set $\Bb{N}_{\fo{\eurm{s}}}$
defined as in  \eqref{1th3} is finite because $D_{r-1} \neq 0$ means that the first $r$ rows and the first  $r$ columns of the matrix $( s_{i + j } )_{i, j = 0}^{\infty}$ are linearly independent. There exists therefore a maximal element $r_{*}:= \max  \Bb{N}_{\fo{\eurm{s}}} \geq 1$ of $\Bb{N}_{\fo{\eurm{s}}}$ for which we have $d_{r_{*}} = \infty$ by virtue of \eqref{1th12}. Then \eqref{1th11} yields  $\eurm{s} =  {\eurm{s}}^{(r_{*})} $  and \eqref{1th5} implies $ \rank ( s_{i + j } )_{i, j = 0}^{\infty} = r_{*}$. Conversely, if $\eurm{s} =  {\eurm{s}}^{(r)} $ for a certain $r \in \Bb{N}_{\fo{\eurm{s}}} $, then  we have the validity of \eqref{1th5} in view of  Theorem~\ref{firsttheorem}\! (b). Thus, $ \rank ( s_{i + j } )_{i, j = 0}^{\infty} = r$ iff $D_{r-1} \neq 0$ and  $\eurm{s} =  {\eurm{s}}^{(r)} $. Combining this assertion with  \eqref{1th7}, \eqref{1th13} and with
\eqref{ioh2}, \eqref{ioh7}, \eqref{6p1th} we obtain   the validity of the following corollary  which contains the statements of Theorems~\ref{thb},~\ref{thc} and Theorem~\ref{thd}\! (a), while Theorem~\ref{thd}\! (b) follows  from \eqref{ioh11} below.

 \vspace{0.2cm}
 \begin{corollary}\label{1cor}
 Let $\eurm{s} := \{s_n \}_{n\geq 0}$ be a  nonzero sequence of real numbers,  $ D_{n} := \det ( s_{i + j } )_{i, j = 0}^{n}$, $n\geq 0$,  and  let the sequence $ {\eurm{s}}^{(m)}:=  \{s_{n}^{(m)}\}_{n \geq 0}$ for every
 $m \geq 1$ satisfying $D_{m-1}\neq 0$  be defined by the expansion% $\vphantom{\frac{A^{A}}{B}}$
 \begin{gather*}
 \frac{Q_{m} (x)}{P_{m} (x) } =  \sum_{k \geq 0 } \ \, \frac{ s_{k}^{(m)}}{x^{k+1}}    \ , \ \ \ \  \  |x| >
 \max \{\ |z| \ | \ P_{m}(z) = 0\ \} \ .
 \end{gather*}

 \medskip
 \indent
 Then $s_{n}^{(m)} = s_{n}$, $0 \leq n \leq 2 m -1$ for every such $m$, and  the infinite Hankel matrix $\mathcal{H}_{\infty} :=   ( s_{i + j } )_{i, j = 0}^{\infty}$  has a finite rank $r\geq 1$ if and only if $ D_{r - 1 } \neq 0$ and one of the following equivalent condition holds:

\vspace{-0.2cm}
\begin{flalign}\label{1cor1}
 &{\rm{{\bf{(a)}}}} \ \ \    \eurm{s} =  {\eurm{s}}^{(r)} \ {\rm{;}}&  &  & \\[0.3cm]\label{1cor2}
  &{\rm{{\bf{(b)}}}} \ \ \  D_{n} = 0 \ ,  \   \ \ n \geq r \ {\rm{;}} & &  & \\[0.3cm] \label{1cor3}
  &{\rm{{\bf{(c)}}}} \ \ \  \sum_{k \geq 0 } \frac{ s_{k}}{x^{k+1}}   = \frac{Q_{r} (x)}{P_{r} (x) }    \ , \ \ \ \  \  |x| >
 \max \{\ |z| \ | \ z \in \Bb{C} \ , \ \ P_{r}(z) = 0\ \} \ {\rm{;}}  & &  &  \\[0.0cm]\label{1cor4}
   &{\rm{{\bf{(d)}}}} \ \ \
  D_{r - 1 }    s_{r + m}\! + \! \sum\limits_{k=0}^{r-1}\ p_{\, r, k}\,  s_{k + m}\! =\!  0\, , \, \   m \geq r \ , \ \,  \mbox{where}\,  \   P_r (x)\! =\! D_{r - 1 } x^{r}\!+\! \sum_{k=0}^{r-1}  p_{\, r, k}\,  x^{k} \, .\hspace{-0.4cm} & &   &
  \end{flalign}

\noindent
In other words, for arbitrary $r \geq 1$ satisfying $D_{r-1}\neq 0$ we have
\begin{gather}\label{1cor5}
    (a) \ \ \Leftrightarrow \ \ (b) \ \ \Leftrightarrow \ \ (c) \ \ \Leftrightarrow \ \ (d) \ \ \Leftrightarrow \ \ \rank\mathcal{H}_{\infty} = r \ ,
\end{gather}

\noindent
while $\rank\mathcal{H}_{\infty} = r $ implies   $D_{r-1}\neq 0$ for arbitrary positive integer $r$.
 \end{corollary}

\vspace{0.1cm}

 \subsection{\texorpdfstring{Positive semidefiniteness}{Positive semidefiniteness}}\label{chs}
The method used by Darboux in deriving formula \cite[(68), p.413]{dar}, now called the Christoffel-Darboux summation formula, can be applied to the polynomials $P_n$. This leads
  to  the following
formula (see \cite[Theorem 4.5, p.23]{chi}, \cite[Theorem 2.6, p.50]{bre})
\begin{gather*}%\label{ioh20}
  D_{r-1}^{2}  \sum_{k=0}^{r-1} \frac{P_{k}(x)P_{k}(y)}{D_{k} D_{k-1}} = \frac{ P_{r}(x)P_{r-1}(y) - P_{r}(y)P_{r-1}(x)}{x-y} \ , \ \ D_{-1} := 1 \ ,  \ x\neq y \ ,
\end{gather*}

\noindent
provided that $D_{k} \neq 0$ for all $0 \leq k \leq r-1$ and $r$ is a positive integer. Thus,
 \begin{gather}\label{ioh21}
  D_{r-1}^{2}  \sum_{k=0}^{r-1} \frac{P_{k}(x)^{2}}{D_{k} D_{k-1}} = P_{r}^{\, \prime }(x)P_{r-1}(x) - P_{r}(x)P_{r-1}^{\, \prime }(x)\ , \
\end{gather}

\noindent
and
\begin{gather}\label{ioh21a}
   D_{r-1}^{2} \sum_{k=0}^{r-1} \frac{\left|P_{k}(z)\right|^{2}}{D_{k} D_{k-1}} =
  \frac{ \im P_{r}(z)P_{r-1}(\overline{z})}{ \im z} \ , \ z \in \Bb{C} \setminus \Bb{R} \ ,
\end{gather}

\noindent
where $\overline{z}$ is a complex conjugate of $z$.

Under the conditions
\begin{gather}\label{ioh14}
 D_{0} > 0 \ , \ \ D_{1} > 0 \ , \ \ldots \ , \   D_{r-1} > 0  \ ,
\end{gather}

\noindent
we then get
\begin{gather*}
 \frac{ \im P_{r}(z)P_{r-1}(\overline{z})}{ \im z}\! =\! \frac{ \im P_{r}(z)\overline{P_{r-1}({z})}}{ \im z}\!=\!
 D_{r-1}^{2} \sum_{k=0}^{r-1} \frac{\left|P_{k}(z)\right|^{2}}{D_{k} D_{k-1}} \!\geq\!  \frac{D_{r-1}}{ D_{r-2}} \left|P_{r-1}(z)\right|^{2} \ , \  \im z \neq 0 \ ,
\end{gather*}

\noindent
so if $P_r(z)=0$ for a $z$ with $\im z\ne 0$, the last inequality above implies $P_{r-1}(z)=0$, which is contradicting that $P_{r}$ and $P_{r-1}$ have no common zeros according to \eqref{int11}.
Therefore all zeros of $P_{r}$  are real, and  if $\lambda$ is a real zero of the polynomial $P_{r}$ we have from \eqref{ioh21}
\begin{gather}\label{ioh22}
   P_{r}^{\, \prime }(\lambda)P_{r-1}(\lambda) = D_{r-1}^{2}  \sum_{k=0}^{r-1} \frac{P_{k}(\lambda)^{2}}{D_{k} D_{k-1}} \ ,
\end{gather}

\noindent
while \eqref{int11} yields
\begin{gather}\label{ioh13}
 P_{r-1} (\lambda) Q_{r} (\lambda) = D_{r-1}^{2} \ .
\end{gather}

\noindent
Since \eqref{ioh13} means  that $P_{r-1} (\lambda) \neq 0$,  \eqref{ioh22} implies  $ P_{r}^{\, \prime }(\lambda) \neq 0$
because
\begin{gather*}
   \sum_{k=0}^{r-1} \frac{P_{k}(\lambda)^{2}}{D_{k} D_{k-1}} \geq \frac{P_{r-1}(\lambda)^{2}}{D_{r-1} D_{r-2}} > 0 \ .
\end{gather*}

\noindent
Thus, all zeros $\{\lambda_{n}\}_{n=1}^{r}$ of $P_{r}$ are simple and by virtue of \eqref{ioh22} and \eqref{ioh13} we have
\begin{gather*}%\label{ioh15}
    \mu_{n}:= \frac{Q_{r} (\lambda_{n})}{P_{r}^{\, \prime }(\lambda_{n})} =
  \left(\sum\limits_{k=0}^{r-1} \dfrac{P_{k}(\lambda_{n})^{2}}{D_{k} D_{k-1}} \right)^{-1} \in (0, +\infty)  \ , \ \ \ \ 1 \leq n \leq r \ ,
\end{gather*}

\noindent
 which gives the following form of the partial fraction decomposition of $Q_{r}/P_{r}$:
\begin{gather}\label{ioh16}
    \frac{Q_{r}(x)}{P_{r}(x)} = \sum_{n=1}^{r} \frac{Q_{r} (\lambda_{n})}{P_{r}^{\, \prime }(\lambda_{n}) (x-\lambda_{n})} =
    \sum_{m=0}^{\infty} \frac{1}{x^{m+1}} \sum_{n=1}^{r} \mu_{n}\lambda_{n}^{m} \ , \ \ \ \  |x| >   \max_{1 \leq n \leq r} \, |\lambda_{n}| \ .
\end{gather}

\medskip Assume now that \eqref{tha1} and \eqref{ioh14} hold. Then $r \in \Bb{N}_{\fo{\eurm{s}}}$ and  \eqref{1cor5} implies the validity of \eqref{1cor3} which in view of \eqref{ioh16} yields
\begin{gather*}%\label{ioh23}
  s_{m} = \sum_{n=1}^{r} \mu_{n}\lambda_{n}^{m} \ , \ \ m \geq 0 \ , \  \mu_{n} = \left(\sum\limits_{k=0}^{r-1} \dfrac{P_{k}(\lambda_{n})^{2}}{D_{k} D_{k-1}} \right)^{-1} > 0 \ ,  \ P_{r} (\lambda_{n}) = 0 \ , \ 1 \leq n \leq r \ .
\end{gather*}

\noindent
We have completely proved the following assertion\footnote{During the preparation of the present paper the second author learned that the result is  formulated
in \cite[Theorem 1.2, p.5]{sh}}
\begin{thmx}[{{\bf{2015}}, \cite[Theorem 1.1, p.1569]{ber1}}]\label{the}
Let $\{s_n \}_{n\geq 0}$ be an arbitrary sequence of real numbers and
 ${\mathcal{H}}_{n} :=   ( s_{i + j } )_{i, j = 0}^{n}$, $n\geq 0$. Assume that there exists a positive integer $n_{0}$ such that

 \vspace{-0.2cm}
\begin{gather*}%\label{fthg1}
 D_n :=   \det {\mathcal{H}}_{n} > 0 \ , \ \ \  0 \leq n \leq n_{0} - 1 \ , \ \  \det {\mathcal{H}}_{n} = 0 \ , \ \ \  n \geq n_{0} \ .
\end{gather*}

\vspace{0.2cm}
\noindent
Then there exist $n_{0}$ distinct real numbers $\{x_{k}\}_{k=1}^{n_{0}}$ and  $n_{0}$ positive numbers
 $\{\mu_{k}\}_{k=1}^{n_{0}}$ such that
 \begin{gather*}%\label{fthg2}
    s_{n}= \int\nolimits_{-\infty}^{+\infty} x^{n} d \mu (x) \ , \ \ \ \ n \geq 0 \ , \ \ \ \
    \mu  := \sum_{k=1}^{n_{0}} \mu_{k} \delta_{x_{k}} \ ,
 \end{gather*}
where $\delta_{y}$ denotes the Dirac measure placed at the point $y \in \Bb{R}$.
\end{thmx}

%%%%%%%%%%%%%%%%%%%%%%%%%%%%%%%%%%%%%%%%%%%%%%%%%%%%%%%%%%%%%%%%%%%%%%%%%%%%%%%%%%%%%%%%%
%%%%%%%%%%%%%%%%%%%%%%%%%%%%%%%%%%%%%%%%%%%%%%%%%%%%%%%%%%%%%%%%%%%%%%%%%%%%%%%%%%%%%%%%%

\vspace{0.35cm}

\section{\texorpdfstring{Proof of Theorem~\ref{firsttheorem}}{Proof of Theorem 3 }}

\vspace{0.25cm}

\subsection{\texorpdfstring{Proof of Theorem~\ref{firsttheorem}(a)}{Proof of Theorem 3(a) } }\label{p1tha} If \eqref{1th1} holds then the first column in the matrices ${\mathcal{H}}_{0}$, ..., ${\mathcal{H}}_{n-1}$ is zero,  and therefore $  D_{0}= D_{1}  =\ \ldots \  =D_{n -1} = 0$, while
\begin{gather}\label{p1tha1}
   D_{n} =
\begin{vmatrix}   0  &   0 &  & \ldots   0 &  s_{n} \\   0  &   0 & \ldots &  s_{n}&  s_{n+1} \\     \ldots  &  \ldots &  \ldots &  \ldots &  \ldots  \\
  0  & s_{n}& \ldots &  s_{2n-2}&  s_{2n-1}  \\    s_{n} & s_{n+1}& \ldots &  s_{2n-1}&  s_{2n}  \\ \end{vmatrix} =
 (-1)^{\tfrac{n (n+1)}{2}} s_{n}^{n+1} \ .
\end{gather}

Conversely, the identity $D_{0}= s_{0}$ together with the condition $D_{0}= 0$  imply $s_{0}=0$. Then \eqref{p1tha1} gives $ D_{1} = - s_{1}^{2}$, which by virtue of the condition  $D_{1}= 0$ yields $s_{1} = 0$.
Pursuing a finite number of repetitions of this fact, we arrive at
  $ s_{0}= s_{1} = \ \ldots \  = s_{n -1} = 0$. In view of \eqref{p1tha1},  $ D_{n} =  (-1)^{\frac{n (n+1)}{2}} s_{n}^{n+1} \neq 0$ and therefore $s_{n} \neq 0$. This concludes the proof of Theorem~\ref{firsttheorem}\! (a).

\vspace{0.15cm}
\subsection{\texorpdfstring{Proof of Theorem~\ref{firsttheorem}(b)}{Proof of Theorem 3(b) }}\label{p1thb}

Since $\eurm{s}$ is a  nonzero sequence, Theorem~\ref{firsttheorem} \!(a) implies that $\Bb{N}_{\fo{\eurm{s}}}$ is nonempty. To prove that the Hankel rank of $ {\eurm{s}}^{(r)}$ is equal to $r$, we observe that the relations \eqref{ioh7} and \eqref{ioh4}
 can also be written as follows
\begin{gather*}%\label{ioh7a}
\begin{pmatrix}
    s_{m +r  }^{(r)} \\    s_{m +r -1 }^{(r)} \\ s_{m +r -2 }^{(r)} \\ \vdots \\
 s_{m +2 }^{(r)} \\ s_{m+1  }^{(r)}
  \end{pmatrix} =
\begin{pmatrix}
 d_{r-1}^{(r)} &d_{r-2}^{(r)}  & \ldots& d_{2}^{(r)} & d_{1}^{(r)} & d_{0}^{(r)} \\[0.1cm]
  1 & 0 &  \ldots & 0 & 0 & 0 \\[0.1cm]
    0 & 1 &  \ldots& 0 & 0 & 0 \\[0.1cm]
   \vdots &  \vdots &  \ldots&  \vdots &  \vdots &  \vdots \\[0.1cm]
  0 & 0 &  \ldots& 1 & 0 & 0 \\[0.1cm]
  0 & 0 &  \ldots& 0 & 1 & 0 \\
\end{pmatrix}
  \begin{pmatrix}
    s_{m +r -1 }^{(r)} \\    s_{m +r -2 }^{(r)} \\ \vdots \\ s_{m +2 }^{(r)} \\
 s_{m +1 }^{(r)} \\ s_{m  }^{(r)}
  \end{pmatrix} \ \ ,
  \ \ \ m \geq 0 \ .
\end{gather*}

\noindent
Therefore
\begin{gather}\label{ioh7b}\begin{pmatrix}
    s_{m +n +r  }^{(r)} \\    s_{m +n +r -1 }^{(r)} \\ s_{m +n +r -2 }^{(r)} \\ \vdots \\
 s_{m +n +2 }^{(r)} \\ s_{m +n +1  }^{(r)}
  \end{pmatrix} =
\begin{pmatrix}
 d_{r-1}^{(r)} &d_{r-2}^{(r)}  & \ldots& d_{2}^{(r)} & d_{1}^{(r)} & d_{0}^{(r)} \\[0.1cm]
  1 & 0 &  \ldots & 0 & 0 & 0 \\[0.1cm]
    0 & 1 &  \ldots& 0 & 0 & 0 \\[0.1cm]
   \vdots &  \vdots &  \ldots&  \vdots &  \vdots &  \vdots \\[0.1cm]
  0 & 0 &  \ldots& 1 & 0 & 0 \\[0.1cm]
  0 & 0 &  \ldots& 0 & 1 & 0 \\
\end{pmatrix}^{n+1}
  \begin{pmatrix}
    s_{m +r -1 }^{(r)} \\    s_{m +r -2 }^{(r)} \\ \vdots \\ s_{m +2 }^{(r)} \\
 s_{m +1 }^{(r)} \\ s_{m  }^{(r)}
  \end{pmatrix} \ \ ,
  \ \ \ m , n \geq 0 \ ,
\end{gather}

\noindent
and the first row of the matrix in the righthand side of \eqref{ioh7b}
gives  the existence of  $r$ numbers  $d_{n, 0}^{(r)}$,  $d_{n, 1}^{(r)}$, ..., $d_{n, r-1}^{(r)}$ satisfying
\begin{gather}\label{ioh11}
  \sum\nolimits_{k=0}^{r-1}\ d_{n, k}^{(r)}\ s_{k + m}^{(r)} = \ s_{r + n +  m}^{(r)} \ ,  \ \ \ \ \ m \geq 0 \ ,
\end{gather}

\noindent
where $d_{0, k}^{(r)}$ is equal to $d_{k}^{(r)}$ from \eqref{ioh7}  for each  $0 \leq k \leq r-1$. This means that
\begin{gather}\label{ioh8}
 d_{n, 0}^{(r)}  \begin{pmatrix}\!
    s_{0}\\ \! s_{1}\\ \! \vdots  \\ \!  s_{r}\\ \! s_{r+1}\\  \! s_{r+2}\\  \! \vdots \\ \!  s_{2r-1}\\  \! s_{2r}^{(r)}\\ \! \! \!   \vdots \\
    \end{pmatrix}\! +\!
d_{n, 1}^{(r)}  \begin{pmatrix}
  \!  s_{1} \! \\ \!   s_{2} \! \\ \!  \vdots \!  \\ \!   s_{r+1} \! \\ \!  s_{r+2} \! \\ \!  s_{r+3} \! \\ \!  \vdots \!  \\ \!  s_{2r}^{(r)} \! \\ \! s_{2r+1}^{(r)} \! \\  \! \vdots \!  \\
    \end{pmatrix}\! + \!...\! +\!
d_{n, r-2}^{(r)}  \begin{pmatrix}\!
    s_{r-2} \! \\ \!   s_{r-1} \! \\ \!  \vdots  \!  \\ \!   s_{2r-2} \! \\ \!  s_{2r-1} \! \\ \! s_{2r}^{(r)} \! \\ \!  \vdots
   \!   \\  \! s_{3r-3}^{(r)} \! \\  \! s_{3r-2}^{(r)} \! \\ \! \vdots \!  \\
    \end{pmatrix}\! +\!
d_{n, r-1}^{(r)}  \begin{pmatrix}\!
    s_{r-1}\\  \!  s_{r} \! \\ \!  \vdots  \!  \\ \!   s_{2r-1} \!   \\ \! s_{2r}^{(r)} \! \\  \! s_{2r+1}^{(r)} \! \\ \! \vdots \!  \\
    \!   s_{3r-2}^{(r)} \! \\ \!  s_{3r-1}^{(r)} \! \\ \! \vdots \!  \\
    \end{pmatrix}\! =\!
 \begin{pmatrix}\!
    s_{r+n}^{(r)}\\ \!   s_{1+r+n}^{(r)} \! \\  \! \vdots  \!  \\  \!  s_{r+r+n}^{(r)} \! \\  \!  s_{r+1+r+n}^{(r)} \! \\  \! s_{r+2+r+n}^{(r)} \! \\ \!    \vdots \!  \\ \!
    s_{2r-1+r+n}^{(r)} \! \\ \! s_{2 r+r+n}^{(r)} \! \\ \!  \vdots  \! \\
    \end{pmatrix} \ ,
\end{gather}

\noindent i.e., for arbitrary $n \geq 0$ the $(r +n+1)$-th column  $( s_{ n + r + j}^{(r)} )_{j= 0}^{\infty}$ of the infinite matrix $( s_{i + j }^{(r)} )_{i, j = 0}^{\infty}$
is a linear combination of the first $r$  columns   which are linearly independent by virtue of \eqref{ioh3}. Thus,
\begin{gather*}%\label{ioh9}
  {\rm{rank}}\,  \big( s_{i + j }^{(r)} \big)_{i, j = 0}^{\infty} = r \ ,
\end{gather*}

\noindent
and we conclude  that the Hankel rank of $\eurm{s}^{(r)}$ is equal to $r$. Theorem~\ref{firsttheorem} \!(b) is proved.

\vspace{0.15cm}
\subsection{\texorpdfstring{Proof of Theorem~\ref{firsttheorem}(d)}{Proof of Theorem 3(d) } }\label{p1thd}

For the sequence $\widehat{\eurm{s}}^{\, (r)} = \{ \widehat{s}_{n}^{\, (r)}\}_{n \geq 0}$ defined by
\begin{gather*}%\label{p1thd1}
\widehat{s}_{n}^{\, (r)} := s_{n} - s_{n}^{(r)} \ , \  n \geq 0 \ , \
  \widehat{\eurm{s}}^{\, (r)} = \eurm{s} -  \eurm{s}^{(r)} \ ,
\end{gather*}

\noindent
we have, by virtue of \eqref{ioh6},
\begin{gather*}%\label{p1thd2}
   \widehat{\eurm{s}}^{\, (r)} = \{ s_{n} - s_{n}^{(r)}\}_{n \geq 0} = \{ \ \underbrace{ 0 \ , \ 0\ ,
   \ 0 \ , \ \ldots \ , \ 0 }_{2r}\ , \
\widehat{s}_{2 r}^{\, (r)}\ , \ \widehat{s}_{2 r+1}^{\, (r)} \ , \  ... \  \}  \ .
\end{gather*}

It is appropriate at this point to recall (see \cite[Definition 8, p.61]{gan}) that two square matrices $A$ and $B$ are called equivalent if there exist two square matrices $P$ and $Q$ with nonzero determinants such that $B = P A Q$. If $\det P = \det Q = 1$ we say that $A$ and $B$ are  {\emph{${\rm{1}}$-equivalent}} and write\vspace{-0.3cm}
\begin{gather*}
  A\  {\overset{1}{\sim}}\ B \ .
\end{gather*}

 \noindent
 It is evident that for two $1$-equivalent square matrices $A$ and $B$ we have $\det A = \det B$,   and also  ${\rm{rank}}\, A = {\rm{rank}}\, B$ in view of \cite[Theorem 2, p.62]{gan}.

\bigskip For arbitrary $r \in \Bb{N}_{\fo{\eurm{s}}}$ consider the matrix
\begin{gather*}
  {\mathcal{H}}_{r}\, =  \! \left(
\begin{array}{llllll}
s_{0} \!\! &\!s_{1} \! \! \! &\!\!\ldots \! & \! s_{r-2}   \!  &\!s_{r-1}   \! &\!s_{r}   \\
s_{1}  \! \! &\! s_{2}\! \! \! &\!\!\ldots \! &\!s_{r-1}   \! &\!s_{r}   \! &\!s_{r+1}   \\
\ldots  \! \! &\! \ldots  \!\! \! &\!\!\ldots \! &\!\ldots \!  &\!\ldots \! &\!\ldots  \\
s_{r-1} \!    \! &\!s_{r} \! \! \! &\!\! \ldots \! &\! s_{2r-3}  \! &\! s_{2r-2}  \! &\! s_{2r-1}  \\
s_{r}   \!  \! &\!s_{r+1} \! \! \! &\!\! \ldots \! &\! s_{2r-2}   \! &\! s_{2r-1}   \! &\! s_{2r}
    \end{array}\!\!
\right)
=  \! \left(
\begin{array}{llllll}
s^{(r)}_{0} \!\! &\!s^{(r)}_{1} \! \! \! &\!\!\ldots \! & \! s^{(r)}_{r-2}   \!  &\!s^{(r)}_{r-1}   \! &\!s^{(r)}_{r}   \\
s^{(r)}_{1}  \! \! &\! s^{(r)}_{2}\! \! \! &\!\!\ldots \! &\!s^{(r)}_{r-1}   \! &\!s^{(r)}_{r}   \! &\!s^{(r)}_{r+1}   \\
\ldots  \! \! &\! \ldots  \!\! \! &\!\!\ldots \! &\!\ldots \!  &\!\ldots \! &\!\ldots  \\
s^{(r)}_{r-1} \!    \! &\!s^{(r)}_{r} \! \! \! &\!\! \ldots \! &\! s^{(r)}_{2r-3}  \! &\! s^{(r)}_{2r-2}  \! &\! s^{(r)}_{2r-1}  \\
s^{(r)}_{r}   \!  \! &\!s^{(r)}_{r+1} \! \! \! &\!\! \ldots \! &\! s^{(r)}_{2r-2}   \! &\! s^{(r)}_{2r-1}   \! &\! s_{2r}
    \end{array}\!\!
\right) \ .
\end{gather*}

\noindent
Subtracting from the last column a linear combination of the first $r$ columns with the coefficients from \eqref{ioh8} with $n=0$, we conclude that
\begin{gather}\label{p1thd3}
   {\mathcal{H}}_{r} \overset{1}{\sim}\begin{pmatrix}
      s_{0} \!\! & \!  s_{1} \!\! & \!  \ldots \!\! & \!  {{s_{r-1}}} \!\! & \!  0  \\
     s_{1}   \!\! & \! s_{2} \!\! & \!   \ldots \!\! & \!  {{s_{r}}} \!\! & \!  0  \\
  \ldots \!\! & \!   \ldots \!\! & \!   \ldots \!\! & \! \ldots \!\! & \!  \ldots \\
{{s_{r-1}}} \!\! & \!  {{s_{r}}} \!\! & \!  \ldots \!\! & \!  {{s_{2r-2}}} \!\! & \! 0   \\
s_{r}   \!\! & \! s_{r+1} \!\! & \!   \ldots \!\! & \!  s_{2r-1} \!\! & \!  \widehat{s}_{2r}^{\, (r)}  \\
    \end{pmatrix}
    \ .
\end{gather}

\noindent
The Laplace expansion of the determinant of the righthand side of  \eqref{p1thd3}  by minors along column $r+1$ leads to the validity of the lefthand equality in \eqref{1th8},
\begin{gather}\label{p1thd4}
D_{r} = D_{r-1} \cdot \left(s_{2r} - s_{2r}^{(r)} \right) \ , \ \ r \in \Bb{N}_{\fo{\eurm{s}}} \ .
\end{gather}

\noindent
Thus, $D_{r} = 0$ if and only if $s_{2 r} = {s}_{2 r}^{(r)}$, which proves \eqref{1th9}  for $d=1$.

To prove \eqref{1th9} for $d >1$ assume  that
\begin{gather}\label{p1thd5}
  s_{2r}\! =\!  s_{2r}^{(r)} \, , \  s_{2r+1}\! =\!  s_{2r+1}^{(r)} \, , \ \ldots \  , s_{2r+d-1}\! = \! s_{2r+d-1}^{(r)} \ .
\end{gather}

\noindent
Consider the matrix
\begin{gather*}
   {\mathcal{H}}_{r +  d} =
 \begin{pmatrix}
      s_{0}^{\, (r)} \!\! & \!  s_{1}^{\, (r)} \!\! & \!  \ldots \!\! & \!  {{s_{r-1}^{\, (r)}}} \!\! & \!  s_{r}^{\, (r)} \!\! & \!  s_{r+1}^{\, (r)}  \!\! & \!  \ldots \!\! & \!  s_{r+d -1 }^{\, (r)} \!\! & \!  s_{r+d  }^{\, (r)}  \\
     s_{1}^{\, (r)}   \!\! & \! s_{2}^{\, (r)} \!\! & \!   \ldots \!\! & \!  {{s_{r}^{\, (r)}}} \!\! & \!  s_{r+1}^{\, (r)} \!\! & \!  s_{r+2}^{\, (r)} \!\! & \!  \ldots \!\! & \!  s_{r+d  }^{\, (r)} \!\! & \!  s_{r+d  +1}^{\, (r)}  \\
  \ldots \!\! & \!   \ldots \!\! & \!   \ldots \!\! & \! \ldots \!\! & \!  \ldots \!\! & \!  \ldots \!\! & \!  \ldots \!\! & \!  \ldots  \!\! & \!  \ldots\\
{{s_{r-1}^{\, (r)}}} \!\! & \!  {{s_{r}^{\, (r)}}} \!\! & \!  \ldots \!\! & \!  {{s_{2r-2}^{\, (r)}}} \!\! & \!  s_{2r-1}^{\, (r)} \!\! & \!
s_{2r}^{\, (r)}  \!\! & \!  \ldots \!\! & \!  s_{2 r+d - 2}^{\, (r)} \!\! & \!  s_{2 r+d-1}^{\, (r)}  \\
s_{r}^{\, (r)}   \!\! & \! s_{r+1}^{\, (r)} \!\! & \!   \ldots \!\! & \!  s_{2r-1}^{\, (r)} \!\! & \!  s_{2r}^{\, (r)} \!\! & \!  s_{2r +1}^{\, (r)} \!\! & \!  \ldots\!\! & \!  s_{2r+d-1}^{\, (r)} \!\! & \!  {{s_{2r+d }}}  \\
 s_{r+1}^{\, (r)} \!\! & \! s_{r+2}^{\, (r)} \!\! & \!  \ldots \!\! & \!  s_{2r}^{\, (r)} \!\! & \!  s_{2r +1}^{\, (r)} \!\! & \!  s_{2r +2}^{\, (r)} \!\! & \!  \ldots\!\! & \!  {{s_{2r+d }}} \!\! & \!  s_{2r+d +1}  \\
  \ldots \!\! & \!   \ldots \!\! & \!   \ldots \!\! & \!  \ldots \!\! & \!  \ldots \!\! & \!  \ldots  \!\! & \!  \ldots\!\! & \!  \ldots  \!\! & \!  \ldots \\
s_{r+d-1}^{\, (r)}   \!\! & \! s_{r+d  }^{\, (r)} \!\! & \!   \ldots \!\! & \!  s_{2r+d-2}^{\, (r)} \!\! & \!  s_{2r+d-1}^{\, (r)} \!\! & \!
 {{s_{2r+d}}} \!\! & \!  \ldots
\!\! & \!  s_{2r+2d-2} \!\! & \!  s_{2r+2d-1}  \\
 s_{r+d  }^{\, (r)}   \!\! & \! s_{r+d+1}^{\, (r)} \!\! & \! \ldots \!\! & \!  s_{2r+d-1}^{\, (r)} \!\! & \!  {{s_{2r+d}}} \!\! & \!  s_{2r+d+1}  \!\! & \!  \ldots
\!\! & \!  s_{2r+2d-1} \!\!     &\! s_{2r+2d} \\
    \end{pmatrix} \ .
\end{gather*}

\vspace{0.25cm}
\noindent For every $0 \leq n \leq d $ we subtract from the $(r+n+1)$-th column  a linear combination of the first $r$ columns with the coefficients from the equality \eqref{ioh8}, and we obtain
 \begin{gather}\label{p1thd6}\hspace{-0.2cm}
  {\mathcal{H}}_{r  + d}\overset{1}{\sim}\begin{pmatrix}
      s_{0}^{\, (r)} \!\! & \!  s_{1}^{\, (r)} \!\! & \!  \ldots \!\! & \!  {{s_{r-1}^{\, (r)}}} \!\! & \!  0 \!\! & \!    0   \!\! & \!  \ldots \!\! & \!   0  \!\! & \!    0  \\
     s_{1}^{\, (r)}   \!\! & \! s_{2}^{\, (r)} \!\! & \!   \ldots \!\! & \!  {{s_{r}^{\, (r)}}} \!\! & \!    0  \!\! & \!   0  \!\! & \!  \ldots \!\! & \!    0  \!\! & \!    0   \\
  \ldots \!\! & \!   \ldots \!\! & \!   \ldots \!\! & \! \ldots \!\! & \!  \ldots \!\! & \!  \ldots \!\! & \!  \ldots \!\! & \!  \ldots  \!\! & \!  \ldots\\
{{s_{r-1}^{\, (r)}}} \!\! & \!  {{s_{r}^{\, (r)}}} \!\! & \!  \ldots \!\! & \!  {{s_{2r-2}^{\, (r)}}} \!\! & \!    0  \!\! & \!   0   \!\! & \!  \ldots \!\! & \!    0  \!\! & \!    0   \\
s_{r}^{\, (r)}   \!\! & \! s_{r+1}^{\, (r)} \!\! & \!   \ldots \!\! & \!  s_{2r-1}^{\, (r)} \!\! & \!   0  \!\! & \!    0  \!\! & \!  \ldots\!\! & \!    0  \!\! & \!  {{\widehat{s}_{2r+d }^{\, (r)}}}  \\
 s_{r+1}^{\, (r)} \!\! & \! s_{r+2}^{\, (r)} \!\! & \!  \ldots \!\! & \!  s_{2r}^{\, (r)} \!\! & \!    0  \!\! & \!    0  \!\! & \!  \ldots\!\! & \!  {{\widehat{s}_{2r+d }^{\, (r)}}} \!\! & \!  \widehat{s}_{2r+d +1}^{\, (r)}  \\
  \ldots \!\! & \!   \ldots \!\! & \!   \ldots \!\! & \!  \ldots \!\! & \!  \ldots \!\! & \!  \ldots  \!\! & \!  \ldots\!\! & \!  \ldots  \!\! & \!  \ldots \\
s_{r+d-1}^{\, (r)}   \!\! & \! s_{r+d  }^{\, (r)} \!\! & \!   \ldots \!\! & \!  s_{2r+d-2}^{\, (r)} \!\! & \!    0  \!\! & \!  {{\widehat{s}_{2r+d}^{\, (r)}}} \!\! & \!  \ldots
\!\! & \!  \widehat{s}_{2r+2d-2}^{\, (r)} \!\! & \!  \widehat{s}_{2r+2d-1}^{\, (r)}  \\
 s_{r+d  }^{\, (r)}   \!\! & \! s_{r+d+1}^{\, (r)} \!\! & \! \ldots \!\! & \!  s_{2r+d-1}^{\, (r)} \!\! & \!  {{\widehat{s}_{2r+d}^{\, (r)}}} \!\! & \!  \widehat{s}_{2r+d+1}^{\, (r)}  \!\! & \!  \ldots
\!\! & \!  \widehat{s}_{2r+2d-1}^{\, (r)} \!\!     &\! \widehat{s}_{2r+2d}^{\, (r)} \\
    \end{pmatrix}  \, .
 \end{gather}

\vspace{0.25cm}
\noindent
Expanding the determinant of the matrix in the righthand side of \eqref{p1thd6} after the last row $d+1$ times or by using that the matrix is
quasi-triangular (see \cite[p.43]{gan}) and using the formulas
\cite[(67), p.43]{gan}, \eqref{p1thd5} and \eqref{ioh2}, we get
\begin{gather}\label{p1thd7}
     D_{r +  d} = (-1)^{{\tfrac{d(d  +1)}{2}} }
\left(s_{2r+d} - s_{2r+d}^{(r)}\right)^{d  +1}  D_{r-1} \  .
\end{gather}

\noindent
 Furthermore, it follows from \eqref{p1thd6} that for all $0 \leq n \leq  d -1$ the matrix ${\mathcal{H}}_{r  + n}$ is 1-equivalent to the matrix with zero $(r+1)$-th column. Therefore
 \begin{gather}\label{p1thd9}
 D_{r} =\  ... \ = D_{r + d-1}=0 \ ,
 \end{gather}
 which proves the implication $\Leftarrow$ in \eqref{1th9}  for $d \geq 1$.

To prove the inverse implication in  \eqref{1th9} for such $d$  assume that \eqref{p1thd9} holds for some $1\leq d < \infty $.
Then $D_{r} = 0$ implies $s_{2r} =  s_{2r}^{(r)}$ by virtue of \eqref{p1thd4}. Therefore \eqref{p1thd5} holds for $d=1$, and we can use the expression \eqref{p1thd7} for $D_{r+1}$   to give $s_{2r+1} =  s_{2r+1}^{(r)}$ if $D_{r+1} = 0$. Pursuing a finite number of repetitions of this trick, we get at last $s_{2r+d-1} =  s_{2r+d-1}^{(r)}$ which completes the proof of \eqref{1th9}.

 Finally, the equivalence  \eqref{1th9} and the implication $\eqref{p1thd5} \Rightarrow \eqref{p1thd7}$ just deduced give the validity of \eqref{1th10}.

To prove the righthand equality in \eqref{1th8} we take  $r \in \Bb{N}_{\fo{\eurm{s}}}$ and consider
\begin{gather*}
  D_{r+1}^{\, \prime}  \! :=  \!\! \left|
\begin{array}{llllll}\!
s_{0} \!\! &\!s_{1} \! \! \! &\!\!\ldots \! & \! s_{r-2}   \!  &\!s_{r-1}   \! &\!s_{r+1} \!  \\ \!
s_{1}  \! \! &\! s_{2}\! \! \! &\!\!\ldots \! &\!s_{r-1}   \! &\!s_{r}   \! &\!s_{r+2}  \! \\ \!
\ldots  \! \! &\! \ldots  \!\! \! &\!\!\ldots \! &\!\ldots \!  &\!\ldots \! &\!\ldots \!\\ \!
s_{r-2} \!    \! &\!s_{r-1} \! \! \! &\!\! \ldots \! &\! s_{2r-4}  \! &\! s_{2r-3}  \! &\! s_{2r-1} \!  \\ \!
s_{r-1} \!    \! &\!s_{r} \! \! \! &\!\! \ldots \! &\! s_{2r-3}  \! &\! s_{2r-2}  \! &\! s_{2r}\!  \\ \!
s_{r}   \!  \! &\!s_{r+1} \! \! \! &\!\! \ldots \! &\! s_{2r-2}   \! &\! s_{2r-1}   \! &\! s_{2r+1}\!
    \end{array}\!\!
\right|
\!\!=  \!\! \left|\begin{array}{llllll}
s^{(r)}_{0} \!\! &\! s^{(r)}_{1} \! \! \!   &\!\!\ldots  \! & \! s^{(r)}_{r-2} \! & \! s^{(r)}_{r-1}    \! & \! s^{(r)}_{r+1}  \!\! \\ \!
s^{(r)}_{1}  \!\!  &\! s^{(r)}_{2} \! \! \!  &\!\!\ldots \! & \! s^{(r)}_{r-1}  \! & \! s^{(r)}_{r}    \! & \! s^{(r)}_{r+2}  \!\! \\ \!
\ldots  \!\!  &\! \ldots  \! \! \! &\!\!\ldots    \! & \! \ldots  \! & \! \ldots  \! & \! \ldots   \\ \!
s^{(r)}_{r-2} \!\!     &\! s^{(r)}_{r-1} \! \! \!   &\!\! \ldots  \! & \! s^{(r)}_{2r-4} \! & \! s^{(r)}_{2r-3}   \! & \!  s^{(r)}_{2r-1}\!\! \\ \!
s^{(r)}_{r-1} \!\!     &\! s^{(r)}_{r} \! \! \!   &\!\! \ldots  \! & \! s^{(r)}_{2r-3} \! & \! s^{(r)}_{2r-2}   \! & \!
 s_{2r}  \!\!\\ \!
s^{(r)}_{r}   \!\!   &\! s^{(r)}_{r+1} \! \! \!   &\!\! \ldots  \! & \! s^{(r)}_{2r-2}  \! & \! s^{(r)}_{2r-1}   \! & \!  s_{2r+1}
    \end{array}\!\!
\right| \, .
\end{gather*}

\noindent
Subtracting from the last column a linear combination of the first $r$ columns with the coefficients from \eqref{ioh8} with $n=0$, we obtain
\begin{gather}\label{p1thd10}  D_{r+1}^{\, \prime}   =
  \left|
\begin{array}{llllll}
s^{(r)}_{0} \!\! &\! s^{(r)}_{1} \! \! \!   &\!\!\ldots  \! & \! s^{(r)}_{r-2} \! & \! s^{(r)}_{r-1}    \! & \!0   \\
s^{(r)}_{1}  \!\!  &\! s^{(r)}_{2} \! \! \!  &\!\!\ldots \! & \! s^{(r)}_{r-1}  \! & \! s^{(r)}_{r}    \! & \! 0   \\
\ldots  \!\!  &\! \ldots  \! \! \! &\!\!\ldots    \! & \! \ldots  \! & \! \ldots  \! & \! \ldots   \\
s^{(r)}_{r-2} \!\!     &\! s^{(r)}_{r-1} \! \! \!   &\!\! \ldots  \! & \! s^{(r)}_{2r-4} \! & \! s^{(r)}_{2r-3}   \! & \!  0 \\
s^{(r)}_{r-1} \!\!     &\! s^{(r)}_{r} \! \! \!   &\!\! \ldots  \! & \! s^{(r)}_{2r-3} \! & \! s^{(r)}_{2r-2}   \! & \!
 s_{2r} -s^{(r)}_{2r} \\
s^{(r)}_{r}   \!\!   &\! s^{(r)}_{r+1} \! \! \!   &\!\! \ldots  \! & \! s^{(r)}_{2r-2}  \! & \! s^{(r)}_{2r-1}   \! & \!  s_{2r+1}-s^{(r)}_{2r+1}
    \end{array}\!\!
\right| \ .
\end{gather}

\noindent
The Laplace expansion of the determinant in the righthand side of  \eqref{p1thd10}  by minors along the column $r+1$ and \eqref{ioh2}  lead to the validity of the righthand equality in \eqref{1th8}.
This finishes the proof of  Theorem~\ref{firsttheorem} \!(d).

%%%%%%%%%%%%%%%%%%%%%%%%%%%%%%%%%%%%%%%%%%%%%%%%%%%%%%%%%%%%%%%%%%%%%%%%%%
%%%%%%%%%%%%%%%%%%%%%%%%%%%%%%%%%%%%%%%%%%%%%%%%%%%%%%%%%%%%%%%%%%%%%%%%%%
%%%%%%%%%%%%%%%%%%%%%%%%%%%%%%%%%%%%%%%%%%%%%%%%%%%%%%%%%%%%%%%%%%%%%%%%%%
%%%%%%%%%%%%%%%%%%%%%%%%%%%%%%%%%%%%%%%%%%%%%%%%%%%%%%%%%%%%%%%%%%%%%%%%%%
\vspace{0.15cm}

\subsection{\texorpdfstring{Proof of Theorem~\ref{firsttheorem}(e)}{Proof of Theorem 3(e) }}\label{p1the}

\vspace{0.15cm} The formula \eqref{p1thd4} proves that
$D_{r} \neq 0$ if and only if $s_{2 r} \neq {s}_{2 r}^{(r)}$. Thus, the definitions of $d_{r}$ given in \eqref{1th11} and in \eqref{1th12} are the same in the case where $d_{r}=0$.

Assume now that the number $d_{r}$ defined by  \eqref{1th11} is finite and  $1\leq d_{r} < \infty$. Then \eqref{p1thd5} holds for $d = d_{r}$ and $s_{2 r + d_{r}} \neq
{s}_{2 r + d_{r}}^{(r)}$ which by \eqref{p1thd7} means that $    D_{r +  d_{r}}\neq 0$. Thus, $d_{r}$ coincides with the infimum in the righthand side of  \eqref{1th12}.

Conversely, if $d_{r}$ is defined by \eqref{1th12} and $1\leq d_{r} < \infty$ then $    D_{r +  d_{r} }\neq 0$ and \eqref{p1thd9} holds for $d = d_{r}$ as well as \eqref{p1thd5} by virtue of \eqref{1th9}. We can therefore apply the formula \eqref{p1thd7} for $d = d_{r}$ to conclude that  $D_{r +  d_{r}}\neq 0$ yields $s_{ 2r + d_{r}} \neq {s}_{2 r + d_{r}}^{(r)}$. This means that $d_{r}$ coincides with the infimum in the righthand side of  \eqref{1th11}.
Thus, two definitions of $d_{r}$ given in \eqref{1th11} and in \eqref{1th12} also coincide when $1\leq d_{r} < \infty$.

Finally, \eqref{1th13} follows directly from \eqref{1th9} applied for every positive integer $d$.  This completes the proof of Theorem~\ref{firsttheorem}\! (e).

%%%%%%%%%%%%%%%%%%%%%%%%%%%%%%%%%%%%%%%%%%%%%%%%%%%%%%%%%%%%%%%%%%%%%%%%%%
%%%%%%%%%%%%%%%%%%%%%%%%%%%%%%%%%%%%%%%%%%%%%%%%%%%%%%%%%%%%%%%%%%%%%%%%%%
%%%%%%%%%%%%%%%%%%%%%%%%%%%%%%%%%%%%%%%%%%%%%%%%%%%%%%%%%%%%%%%%%%%%%%%%%%
%%%%%%%%%%%%%%%%%%%%%%%%%%%%%%%%%%%%%%%%%%%%%%%%%%%%%%%%%%%%%%%%%%%%%%%%%%

\vspace{0.15cm}
\subsection{\texorpdfstring{Proof of Theorem~\ref{firsttheorem}(c)}{Proof of Theorem 3(c) }}\label{p1thc}

\vspace{0.15cm}
 If $p (x, t) = \sum_{k=0}^{n}\sum_{j=0}^{m} p_{\,k,\, j}\,x^{k} t^{j}$ is an algebraic polynomial with real coefficients of two variables $x$ and $t$,  we use the linear functional $L$ from \eqref{eq:L} with respect to the $t$-variable to get
\begin{gather*}    L_t \big( p(x, t)\big)  = \sum\nolimits_{k=0}^{n}\sum\nolimits_{j=0}^{m} s_{j} p_{k,\, j}\, x^{k} =
 \sum\nolimits_{k=0}^{n} \left(\sum\nolimits_{j=0}^{m} s_{j} p_{\,k,\, j} \right) x^{k}  \in \eusm{P}[\Bb{R}]  \ .
\end{gather*}

\noindent
For example,
\begin{gather*}
   L_t \left(\! \frac{1\! -\!1}{x\!-\!t}\!\right)\! =\! 0 \, ,  \
   L_t  \left(\!\frac{x\! -\!t}{x\!-\!t}\!\right) \! =\! s_{0} \, ,  \
   L_t  \left(\!\frac{x^{2}\! -\!t^{2}}{x\!-\!t}\!\right) \! =\! s_{0}x\! +\! s_{1} \, ,  \
     L_t  \left(\!\frac{x^{n}\! -\!t^{n}}{x\!-\!t}\!\right) \! = \!\sum\limits_{j=0}^{n\!-\!1} s_{j} x^{n\!-\!1\!-\!j} \, , \ n \!\geq\! 1  \, .
\end{gather*}

\noindent
Comparing these equalities with \eqref{thb1} we see that
\begin{gather*}%\label{1p1th}
    Q_{n} (x) := L_t \left(\frac{P_{n} (x) - P_{n} (t)}{x-t}\right) \ , \ \ \ \   n \geq 0 \ ,
\end{gather*}

\noindent
and if  $P_n (x)\! =:\! \sum_{k=0}^{n}  p_{n, k} x^{k}$, $Q_{n} (x)\! =:\! \sum_{k=0}^{n-1}  q_{n, k} x^{k}$, $n \geq 1$,  we obtain
\begin{align*} %\label{2p1th}
 & P_{0}(x)\! =\! 1  \ , \ P_{1}(x) \! =\! s_{0} x \! -\! s_{1} \ , \ \
  P_n (x)\! = \! D_{n-1} x^{n} \!  + \! \sum\limits_{k=0}^{n-1}
  p_{n, k}\, x^{k}   \, , \ p_{n, n}\! = \!D_{n-1}  \, ,   \\ &
\frac{P_{n}(x)\!-\!P_{n}(t)}{x\!-\!t}\! =\!\sum_{m=1}^{n}  p_{n, m} \frac{x^{m}\!-\!t^{m}}{x\!-\!t}  \!=\!
\sum_{m=0}^{n-1}  p_{n, m+1}\sum_{k=0}^{m} x^{k} t^{m-k}\! = \!
\sum_{k=0}^{n-1} \left[\sum_{m=k}^{n-1} p_{n, m+1} t^{m-k}\right]  x^{k}\ ,
 \\  &   Q_{1}(x) := s_{0}^{2} \ , \
Q_{n} (x) =
L_t \left(\frac{P_{n} (x) - P_{n} (t)}{x-t}\right) = \sum_{k=0}^{n-1}\Big( \sum_{m=k}^{n-1} p_{n, m+1}s_{m-k}\Big) x^{k} \, ,
\\ &
   q_{n, m} = \sum_{k=m}^{n-1} p_{n, k+1}s_{k-m} = \sum_{k=0}^{n-m-1} p_{n, k+m+1}s_{k} \ , \  \ \    0 \leq m \leq n-1 \ ,   \ \ \  n \geq 1 \ ,
\end{align*}

\noindent
where the latter equalities can be written in the following form:
\begin{gather}\label{2ap1th}
 q_{n, n-1- m } = \sum_{k=0}^{m} p_{n, n - (m-k) } s_{k} \ ,  \ \ \   0 \leq m \leq n-1 \ ,  \ \ \   n \geq 1 \ .
\end{gather}

\medskip
 Let $r \in \Bb{N}_{\fo{\eurm{s}}}$, i.e., $D_{r-1} \neq 0 $.
According to \eqref{thb1} we have $L ( t^{m} P_r (t) )=0$ for every $0 \leq m \leq r-1$ and therefore
\begin{gather}\label{3p1th}
 \sum_{k=0}^{r}  p_{r, k} s_{k +m} = 0 \ , \  \ \ \   0 \leq m \leq r-1 \ ,
\end{gather}

\noindent
which gives
\begin{gather*}%\label{4p1th}
   \sum_{k=0}^{r-1} \left(- p_{r, k}\right) s_{k +m} =  D_{r-1}  s_{r +m}\ , \  \ \ \   0 \leq m \leq r-1 \ ,
\end{gather*}

\noindent
or,
\begin{gather}\label{5p1th}
   \begin{pmatrix} s_{0} &  s_{1}& \ldots &  s_{r-2} &  s_{r-1} \\[0.1cm]
     s_{1} & s_{2}& \ldots &  s_{r-1}&  s_{r} \\
          \ldots  &  \ldots &  \ldots &  \ldots &  \ldots  \\[0.1cm]
          s_{r-2} & s_{r-1}& \ldots &  s_{2r-4}&  s_{2r-3}  \\[0.1cm]
              s_{r-1} & s_{r}& \ldots &  s_{2r-3}&  s_{2r-2}
              \\ \end{pmatrix}
\begin{pmatrix}
- p_{r, 0}/D_{r-1} \\ - p_{r, 1}/D_{r-1} \\
  \ldots \\ - p_{r, r-2}/D_{r-1} \\  - p_{r, r-1}/D_{r-1} \\
\end{pmatrix} = \begin{pmatrix}
 s_{r} \\[0.1cm]  s_{r+1} \\
  \ldots \\[0.1cm] s_{2r-2} \\[0.1cm]  s_{2r-1} \\
\end{pmatrix} \ .
\end{gather}

\noindent
Since $D_{r-1} \neq 0 $ it follows from \eqref{ioh4} and \eqref{5p1th}  that
\begin{gather}\label{6p1th}
 d_{k}^{(r)} =  - \frac{p_{r, k}}{D_{r-1}}   \ , \  \ \ \   0 \leq k \leq r-1 \ , \  \
  P_r (x)\! = \! D_{r-1} x^{r} \!  + \! \sum\limits_{k=0}^{r-1}
  p_{r, k}\, x^{k}   \ , \ \  r \in \Bb{N}_{\fo{\eurm{s}}} \ ,
\end{gather}

\noindent
and therefore the recursive formulas \eqref{ioh5} can be written in the following manner
\begin{gather*}%\label{6ap1th}
D_{r-1} s_{2 r + m }^{(r)} \!+ \!\sum\limits_{k=0}^{r-1} \ p_{r, k}  s_{r + m + k}^{(r)} = 0 \, ,  \  \ \ \     m \geq 0  \, ,
\end{gather*}

\noindent
which together with \eqref{3p1th} gives
\begin{gather*}%\label{6bp1th}
\sum\limits_{k=0}^{r} \ p_{r, k}  s_{ k + m }^{(r)} = 0 \, ,  \  \ \ \     m \geq 0  \, ,
\end{gather*}

\noindent
or
\begin{gather}\label{6cp1th}
D_{r-1}  s_{ r + m }^{(r)} + \sum\limits_{k=0}^{r-1} \ p_{r, k}  s_{ k + m }^{(r)} = 0 \, ,  \   \ \ \    m \geq 0  \, .
\end{gather}

\medskip
Denote $P_{r}^{\star} (x) := x^{r}P_{r}(1/x)$ and $Q_{r}^{\star} (x) := x^{r-1} Q_{r}(1/x)$.  Then
\begin{gather*}%\label{7p1th}
  P_{r}^{\star} (x)\! = \! D_{r-1}  \!  + \! \sum\limits_{k=0}^{r-1}
  p_{r, k}\, x^{r - k} \ , \ \   P_{r}^{\star} (0)\! = \! D_{r-1}\neq 0 \ , \ \
 Q_{r}^{\star} (x)\! =\! \sum_{k=0}^{r-1}  q_{r, k} x^{r-1- k}   \ ,
\end{gather*}

\noindent
and the function $\psi_{r} (z):=  Q_{r}^{\star} (z)/ P_{r}^{\star} (z)$  is analytic on the open disk $|z|< 1/\rho_{r}$ where $\rho_{r} := \max \{\ |z| \ | \ z \in \Bb{C} \ , \ P_{r}(z) = 0\}$. Let $a_{n}$, $n \geq 0$, be the coefficients of the Taylor expansion of $\psi_{r} (z)$ at the origin,
\begin{gather}\label{9p1th}
  \psi_{r} (z) = \sum_{n \geq 0} a_{n}z^{n} \ , \  \ \ \   |z| < 1/\rho_{r} \ ,
\end{gather}

\noindent
which is obviously equivalent to the expansion of the form
\begin{gather}\label{10p1th}
 \frac{Q_{r} (x)}{P_{r} (x) }   =  \sum_{m \geq 0 } \frac{ a_{m}}{x^{m+1}}  \ , \  \ \ \  |x| > \rho_{r} \ .
\end{gather}

\noindent
The identity
\begin{gather*}%\label{11p1th}
  \left( \sum_{m=0}^{r} p_{r, m} x^{m}\right) \sum_{m \geq 0} \frac{a_{m}}{x^{m+1}} =\sum_{m \geq 0 } \frac{
  \begin{displaystyle}
  \sum\nolimits_{k=0}^{r}
  \end{displaystyle}
   \ p_{r, k} a_{m+k}}{x^{m+1}}    +
\sum_{m=0}^{r-1}x^{m} \sum_{k=m}^{r-1}  p_{r, k+1}  a_{k- m}  \ ,
\end{gather*}

\noindent
and \eqref{10p1th} imply
\begin{equation*}%\label{12p1th}
\begin{split} &
  \sum_{k=m}^{r-1}  p_{r, k+1}  a_{k- m} = q_{r, m} \ , \ \ \ \   0 \leq m \leq r-1 \ , \\ &
    \sum\nolimits_{k=0}^{r} \ p_{r, k} a_{m+k} = 0 \ , \ \ \ \  \ \  m \geq 0 \ ,
\end{split}
\end{equation*}

\noindent
which can be written as
\begin{equation}\label{13p1th}
\begin{split} &
  \sum_{k=0}^{m}  p_{r, r- (m- k)}  a_{k} = q_{r,  r- 1 - m} \ , \  \ \ \  \ \ \   \ \ \   \ \ \    0 \leq m \leq r-1 \ , \\ &
  D_{r-1} a_{m+r} + \sum\nolimits_{k=0}^{r-1} \ p_{r, k} a_{m+k} = 0 \ , \  \ \ \  \ \  m \geq 0 \ .
\end{split}
\end{equation}

\noindent
Observe that \eqref{ioh2}, \eqref{2ap1th} for $n=r$ and \eqref{6cp1th} mean that
\begin{equation*}%\label{14p1th}
\begin{split}&
  \sum_{k=0}^{m}  p_{r, r- (m- k)}  s_{k}^{(r)} = q_{r, r-1- m} \ , \  \ \ \  \ \ \   0 \leq m \leq r-1 \ , \\ &
  D_{r-1}s_{m+r}^{(r)}  + \sum\nolimits_{k=0}^{r-1} \ p_{r, k} s_{m+k}^{(r)} = 0 \ , \ \  m \geq 0 \ .
\end{split}
\end{equation*}

\noindent
and by subtracting from these equalities the corresponding equalities in \eqref{13p1th} we obtain
\begin{equation*}%\label{15p1th}
\begin{split}&  D_{r-1} \left[ s_{m}^{(r)} - a_{m}\right] +
  \sum_{k=0}^{m-1}  p_{r, r- (m- k)}  \left[ s_{k}^{(r)} - a_{k}\right] = 0 \ , \  \ \ \   \ \ \  0 \leq m \leq r-1 \ , \\ &
  D_{r-1}\left[s_{m+r}^{(r)} - a_{m+r}\right]  + \sum\nolimits_{k=0}^{r-1} \ p_{r, k} \left[s_{m+k}^{(r)} - a_{m+k}\right] = 0 \ , \ m \geq 0 \ ,
\end{split}
\end{equation*}

\noindent
where it is assumed that $\sum_{0}^{-1} := 0$. From these recurrence relations we obtain $s_{m}^{(r)} = a_{m}$ for all $m \geq 0$ as a consequence of $ D_{r-1} \neq 0$. Together with  \eqref{ioh2} this proves \eqref{1th7}.
Since the radius  of  convergence of the Taylor series \eqref{9p1th} is  known we get the validity of \eqref{1th6} by virtue of the Cauchy-Hadamard formula (see \cite[(2), p.200]{rud}).  Theorem~\ref{firsttheorem}\! (c) is proved.
%%%%%%%%%%%%%%%%%%%%%%%%%%%%%%%%%%%%%%%%%%%%%%%%%%%%%%%%%%%%%%%%%%%%%%%%%%
%%%%%%%%%%%%%%%%%%%%%%%%%%%%%%%%%%%%%%%%%%%%%%%%%%%%%%%%%%%%%%%%%%%%%%%%%%
%%%%%%%%%%%%%%%%%%%%%%%%%%%%%%%%%%%%%%%%%%%%%%%%%%%%%%%%%%%%%%%%%%%%%%%%%%
%%%%%%%%%%%%%%%%%%%%%%%%%%%%%%%%%%%%%%%%%%%%%%%%%%%%%%%%%%%%%%%%%%%%%%%%%%

\vspace{0.25cm}

\section{\texorpdfstring{Proof of Theorem~\ref{maintheorem}}{Proof of Theorem 1 } }\label{p0th}

\vspace{0.15cm}
If $t_{n}=0$ for all $n \geq 0$ then Theorem~\ref{firsttheorem}\! (a) yields $s_{n}=0$ for every $n \geq 0$ and therefore in this case \eqref{intr1} has only one solution as stated in the theorem.

To examine the case $   Z_{\hspace{0.025cm} \eurm{t}}  \neq \emptyset$ we introduce the notation
\begin{gather}\label{1p0th}
\Delta_{0} \! := \!  (-1)^{\tfrac{n_{0}\! +\! 1}{2}} t_{n_{0}}  \, , \
\Delta_{k+1}\!  :=\!  (-1)^{\tfrac{n_{k\! +\! 1}\! -\! n_{k}}{2}} t_{n_{k\! +\! 1}}  t_{n_{k}}  \, ,  \ 0 \!\leq\! k \!<\! m\! -\!1 \, , \  2\! \leq\! m \!\leq \!\infty \ .
\end{gather}

\vspace{0.15cm}
\subsection{\texorpdfstring{ Necessity of Theorem~\ref{maintheorem}}{ Necessity of Theorem 1 } }

To prove  necessity, we assume that \eqref{intr1} has at least one solution $\eurm{s}\!:= \!\{s_n \}_{n\geq 0}$ for a given $\{t_{n}\}_{n \geq 0}$.
Then by Theorem~\ref{firsttheorem}\! (a),
\begin{gather}\label{2p0th}
t_{n_{0}} = (-1)^{\frac{n_{0} (n_{0}+1)}{2}} s_{n_{0}}^{n_{0}+1} \ .
\end{gather}

\noindent
Furthermore, in \eqref{1th3} we have
\begin{gather*}%\label{5p0th}
   \Bb{N}_{\fo{\eurm{s}}} = \left\{ n_{k}+1\right\}_{0 \leq k < m} \ ,
\end{gather*}

\noindent
and for every $0 \leq k < m$ the formulas \eqref{ioh2}, \eqref{ioh4} and \eqref{ioh5} for $r =  n_{k}+1$
and the numbers $s_{0}$,  $s_{1}$, ..., $s_{2n_{k} +1}$
produce the sequence
\begin{gather}\label{6p0th}
{\eurm{s}}^{(n_{k}+1)} =\{ s_{0} \ , \ s_{1}\ ,  \ s_{2} \ , \ \ldots \ , \ s_{2 n_{k}+1} \ , \
s_{2n_{k}+2}^{(n_{k}+1)}\ , \ s_{2n_{k}+3}^{(n_{k}+1)} \ , \  ... \  \} \ , \ \ 0 \leq k < m \ .
\end{gather}

\noindent For arbitrary $0 \leq k < m -1$,  $2 \leq m \leq \infty$,
it follows from Theorem~\ref{firsttheorem}\! (d) with  $r =n_{k} +1$ and $d= n_{k+1} - n_{k}-1 $
that
\begin{gather}\label{7p0th}
    t_{n_{k+1}} = (-1)^{{\tfrac{(n_{k+1} - n_{k}-1)(n_{k+1} - n_{k})}{2}} } \left(s_{n_{k+1} + n_{k}+1} - s_{n_{k+1} + n_{k}+1}^{(n_{k} +1)}\right)^{n_{k+1} - n_{k}}  t_{n_{k}} \ .
\end{gather}

\noindent
But in view of \eqref{1p0th}, \eqref{2p0th} and \eqref{7p0th} we have
\begin{align*} &
\Delta_{0}  = (-1)^{\frac{ (n_{0}+1)^{2}}{2}} s_{n_{0}}^{n_{0}+1} > 0  \ , \\ &
 \Delta_{k+1}  =
   (-1)^{{\tfrac{(n_{k+1} - n_{k})^{2}}{2}} } \left(s_{n_{k+1} + n_{k}+1} - s_{n_{k+1} + n_{k}+1}^{(n_{k} +1)}\right)^{n_{k+1} - n_{k}}   t_{n_{k}}^{2} > 0 \ , \ 0 \leq k < m-1 \ ,
\end{align*}

\noindent
provided that $n_{0}+1 \in 2 \Bb{N}$ and $ n_{k+1}-n_{k} \in 2 \Bb{N}$ for $0 \leq k < m -1$,  $2 \leq m \leq \infty$. This proves the necessity part of Theorem~\ref{maintheorem}.

\vspace{0.15cm}

\subsection{\texorpdfstring{  Sufficiency of Theorem~\ref{maintheorem}}{  Sufficiency of Theorem 1 }}

The proof of sufficiency proceeds by induction on $k$. Given a sequence $\{t_{n}\}_{n \geq 0}$ satisfying the conditions
of Theorem~\ref{maintheorem} we will determine the terms of the sequence $\eurm{s}\!:= \!\{s_n \}_{n\geq 0}$ such that
\eqref{intr1} holds. Let  $D_{n} := \det (s_{i + j})_{i, j = 0}^{n}$, $n \geq 0$.

If $n_{0}=0$ we put $s_{0} = t_{0}$ to obtain $D_{0} = s_{0}= t_{0}$. If $n_{0}\geq 1$ and  $n_{0}+1 \in 2 \Bb{N} +1$,  we set
 \begin{gather*}%\label{10p0th}
    s_{0}= s_{1}  =\ \ldots \  =s_{n_{0} -1} = 0  \ , \ s_{n_{0}} =  (-1)^{\frac{n_{0} }{2}}  t_{n_{0}}^{\tfrac{1}{n_{0}+1}}
    \  .
\end{gather*}
According to Theorem~\ref{firsttheorem}\! (a) we have
 \begin{gather*}%\label{11p0th}
    D_{0}= D_{1}  =\ \ldots \  =D_{n_{0} -1} = 0  \ , \ D_{n_{0}} = (-1)^{\frac{n_{0} (n_{0}+1)}{2}} s_{n_{0}}^{n_{0}+1} =  (-1)^{n_{0} (n_{0}+1)} t_{n_{0}} =    t_{n_{0}}  \ .
\end{gather*}

\noindent Assume now that  $n_{0}+1 \in 2 \Bb{N}$. Then, in view of \eqref{1p0th},
\begin{gather*}
t_{n_{0}}   =  (-1)^{\tfrac{n_{0}+1}{2}} \Delta_{0} \ , \ \Delta_{0} > 0 \ ,
\end{gather*}

 \noindent
 and if we put
\begin{gather*}%\label{12p0th}
    s_{0}= s_{1}  =\ \ldots \  =s_{n_{0} -1} = 0  \ , \ s_{n_{0}} =   \Delta_{0}^{\tfrac{1}{n_{0}+1}}
    \  ,
\end{gather*}

 \noindent
 then by  Theorem~\ref{firsttheorem}\! (a) we obtain
 \begin{gather*} D_{0}= D_{1}  =\ \ldots \ = D_{n_{0} -1} = 0  \ , \
  D_{n_{0}} =  (-1)^{\tfrac{n_{0} (n_{0}+1)}{2}} s_{n_{0}}^{n_{0}+1} = (-1)^{\tfrac{n_{0} (n_{0}+1)}{2}}\Delta_{0} = t_{n_{0}} \ .
 \end{gather*}

\noindent
Therefore in both cases the numbers  $s_{0}$, ... , $s_{2n_{0}}$ with $s_{n_{0} +1}$, ... , $s_{2n_{0}}$  chosen arbitrarily satisfy \eqref{intr1} for $0 \leq n \leq n_{0}$.

 \medskip
 Suppose that for a certain $k$ satisfying $0 \leq k < m -1$ where  $2 \leq m \leq \infty$, the numbers  $s_{0}$, ... , $s_{2n_{k}}$ satisfy \eqref{intr1} for $0 \leq n \leq n_{k}$.
 We prove that it is possible to determine the numbers  $s_{2n_{k} +1}$,  $s_{2n_{k} +2}$ ... , $s_{2n_{k+1}}$ such that
 \eqref{intr1} holds for $0 \leq n \leq n_{k+1}$.

Choosing arbitrarily the number  $s_{2n_{k} +1}$ we construct the sequence ${\eurm{s}}^{(n_{k}+1)}$ as in \eqref{6p0th}.

 \noindent
 Assume first that $n_{k+1} = n_{k} +1$. In view of \eqref{1th8} for $r = n_{k} +1$,
 \begin{gather*}
    D_{n_{k+1} } = \left(s_{2n_{k+1}} - s_{2n_{k+1}}^{(n_{k}+1)}\right)  t_{n_{k}} \ ,
 \end{gather*}

\noindent
and by putting
\begin{gather*}
  s_{2n_{k+1}} = s_{2n_{k+1}}^{(n_{k}+1)} + \frac{t_{n_{k+1}}}{t_{n_{k}}} \ ,
\end{gather*}

 \noindent
 we obtain the desired equality $D_{n_{k+1} } = t_{n_{k+1}}$.

\medskip
 Assume now that $n_{k+1} - n_{k} \geq 2$. Then
 \begin{gather*}
  t_{n_{k} +1} = ... = t_{ n_{k+1} -1} = 0
 \end{gather*}

 \noindent
and if  we set
  \begin{gather*}%\label{9p0th}
s_{2n_{k} +2}=  s_{2n_{k} +2}^{(n_{k} +1)} \, , \  s_{2n_{k} +3}=  s_{2n_{k} +3}^{(n_{k} +1)} \, , \ \ldots \ , \
s_{n_{k+1} + n_{k} }=s_{n_{k+1} + n_{k}}^{(n_{k} +1)} \ ,
  \end{gather*}

\noindent
we obtain, by virtue of \eqref{1th9} with $r =n_{k} +1$ and $d= n_{k+1} - n_{k}-1 $,
\begin{gather*}
    D_{n_{k} +1} = 0 =  t_{n_{k} +1} \ ,   ... \ ,  D_{ n_{k+1} -1} = 0  = t_{ n_{k+1} -1} \ ,
\end{gather*}

\noindent
and in view of \eqref{1th10},
\begin{gather}\label{31p0th}
    D_{n_{k+1}} = (-1)^{{\tfrac{(n_{k+1} - n_{k}-1)(n_{k+1} - n_{k})}{2}} } \left(s_{n_{k+1} + n_{k}+1} - s_{n_{k+1} + n_{k}+1}^{(n_{k} +1)}\right)^{n_{k+1} - n_{k}}  t_{n_{k}} \ .
\end{gather}

\medskip
If $n_{k+1} - n_{k}\in 2 \Bb{N} +1$ we can choose
\begin{gather*}
   s_{n_{k+1} + n_{k}+1} = s_{n_{k+1} + n_{k}+1}^{(n_{k} +1)} + (-1)^{{\tfrac{(n_{k+1} - n_{k}-1)}{2}} }
  \left( \frac{ t_{n_{k+1}} }{ t_{n_{k}} }\right)^{\tfrac{1}{n_{k+1} - n_{k}}} \ ,
\end{gather*}

\noindent
to have from \eqref{31p0th}, $  D_{n_{k+1}} =   t_{n_{k+1}}$.

\medskip
 But if  $n_{k+1} - n_{k}\in 2 \Bb{N} $ we set
\begin{gather*}
   s_{n_{k+1} + n_{k}+1} = s_{n_{k+1} + n_{k}+1}^{(n_{k} +1)} +
  \left(\frac{\Delta_{k+1}}{t_{n_{k}}^{2}}\right)^{\tfrac{1}{n_{k+1} - n_{k}}} \ ,
\end{gather*}

\noindent
where according to the conditions of the theorem
\begin{gather*}
 \Delta_{k+1}  =  (-1)^{\tfrac{n_{k+1}-n_{k}}{2}} t_{n_{k+1}}  t_{n_{k}} \ , \ \ \Delta_{k+1} > 0 \ .
\end{gather*}

\noindent
Then \eqref{31p0th} gives
 \begin{align*}
    D_{n_{k+1}} & = (-1)^{{\tfrac{(n_{k+1} - n_{k}-1)(n_{k+1} - n_{k})}{2}} } \left(s_{n_{k+1} + n_{k}+1} - s_{n_{k+1} + n_{k}+1}^{(n_{k} +1)}\right)^{n_{k+1} - n_{k}}  t_{n_{k}}  \\  &  =
 (-1)^{{\tfrac{(n_{k+1} - n_{k}-1)(n_{k+1} - n_{k})}{2}} }\frac{\Delta_{k+1}}{t_{n_{k}} }  \\  &  =
 (-1)^{{\tfrac{(n_{k+1} - n_{k}-1)(n_{k+1} - n_{k})}{2}} }  (-1)^{\tfrac{n_{k+1}-n_{k}}{2}} t_{n_{k+1}} =
  (-1)^{\tfrac{(n_{k+1}-n_{k})^{2}}{2}} t_{n_{k+1}} = t_{n_{k+1}} \ .
   \end{align*}

  \noindent
 Therefore in both cases the numbers  $s_{0}$, ... , $s_{n_{k+1} + n_{k}+1}$ satisfy \eqref{intr1} for $0 \leq n \leq n_{k+1}$
 independently on the choice of $s_{n_{k+1} + n_{k}+2}$, ... , $s_{2 n_{k+1}}$. Choosing arbitrarily the latter numbers  we obtain the desired result. This finishes the proof of Theorem~\ref{maintheorem}.

%%%%%%%%%%%%%%%%%%%%%%%%%%%%%%%%%%%%%%%%%%%%%%%%%%%%%%%%%%%%%%%%%%%%%%%%%%%%%%%%%%%%%%%%%
%%%%%%%%%%%%%%%%%%%%%%%%%%%%%%%%%%%%%%%%%%%%%%%%%%%%%%%%%%%%%%%%%%%%%%%%%%%%%%%%%%%%%%%%%

%%%%%%%%%%%%%%%%%%%%%%%%%%%%%%%%%%%%%%%%%%%%%%%%%%%%%%%%%%%%%%%%%%%%%%%%%%%%%%%%%%%%%%%%%%%%%%%%%%%%%%%%
%%%%%%%%%%%%%%%%%%%%%%%%%%%%%%%%%%%%%%%%%%%%%%%%%%%%%%%%%%%%%%%%%%%%%%%%%%%%%%%%%%%%%%%%%%%%%%%%%%%%%%%%
%%%%%%%%%%%%%%%%%%%%%%%%%%%%%%%%%%%%%%%%%%%%%%%%%%%%%%%%%%%%%%%%%%%%%%%%%%%%%%%%%%%%%%%%%%%%%%%%%%%%%%%%
%%%%%%%%%%%%%%%%%%%%%%%%%%%%%%%%%%%%%%%%%%%%%%%%%%%%%%%%%%%%%%%%%%%%%%%%%%%%%%%%%%%%%%%%%%%%%%%%%%%%%%%%

\end{document}